\documentclass[12pt,noamsfonts]{amsart}

\usepackage{amsfonts}
\usepackage{amssymb}   
\usepackage{amscd}

\def\ZZ{{\mathbb Z}}

\def\AA{{\mathbb A}}
\def\PP{{\mathbb P}}
\def\RR{{\mathbb R}}

\def\aa{{\mathbf a}}
\def\bb{{\mathbf b}}
\def\cc{{\mathbf c}}
\def\ee{{\mathbf e}}
\def\ep{{\boldsymbol\varepsilon}}
\def\hh{{\mathbf h}}
\def\pp{{\mathbf p}}
\def\qq{{\mathbf q}}
\def\uu{{\mathbf u}}
\def\vv{{\mathbf v}}
\def\ww{{\mathbf w}}
\def\b00{\bar{\mathbf{0}}}
\def\baa{\bar{\mathbf{a}}}
\def\bbb{\bar{\mathbf{b}}}

\def\bee{\bar{\mathbf{e}}}
\def\buu{\bar{\mathbf{u}}}

\def\cD{\mathcal{D}}
\def\cF{\mathcal{F}}
\def\cG{\mathcal{G}}

\def\cO{\mathcal{O}}

\def\tF{\widetilde{F}}
\def\tG{\widetilde{G}}
\def\Dl{{\it D}_L}
\def\Dr{{\it D}_R}
\def\D{{\it D}}
\def\Ml{\text{{\sf $A$-GrMod$_{\theta}$}}}
\def\Mr{\text{{\sf GrMod$_{\theta}$-$A$}}}
\def\Mlf{\text{{\sf $A$-GrMod$_{\theta}^{f}$}}}
\def\Mrf{\text{{\sf GrMod$_{\theta}^{f}$-$A$}}}
\def\DMod{{\text {\sf $\mathcal{D}$-Mod}}}
\def\ModD{{\text {\sf Mod-$\mathcal{D}$}}}
\def\DCoh{{\text {\sf $\mathcal{D}$-Coh}}}
\def\CohD{{\text {\sf Coh-$\mathcal{D}$}}}
\def\btors{{\text {\sf $\mathfrak{b}$-Tors}}}
\def\torsb{{\text {\sf Tors-$\mathfrak{b}$}}}
\def\hsigma{\widehat{\sigma}}

\DeclareMathOperator{\Ker}{Ker}
\DeclareMathOperator{\Pic}{Pic}
\DeclareMathOperator{\Hom}{Hom}
\DeclareMathOperator{\End}{End}
\DeclareMathOperator{\Spec}{Spec}
\DeclareMathOperator{\Cl}{Cl}
\DeclareMathOperator{\Stab}{Stab}
\DeclareMathOperator{\CharVar}{Ch}
\DeclareMathOperator{\Jac}{Jac}
\DeclareMathOperator{\IN}{in}
\DeclareMathOperator{\gr}{gr}
\DeclareMathOperator{\Ann}{Ann}
\DeclareMathOperator{\Lie}{Lie}

\DeclareMathOperator{\Var}{Var}

\newtheorem{lemma}{Lemma}[section]
\newtheorem{theorem}[lemma]{Theorem}
\newtheorem*{theorem*}{Theorem}
\newtheorem{corollary}[lemma]{Corollary}
\newtheorem{proposition}[lemma]{Proposition}

\theoremstyle{definition}
\newtheorem{example}[lemma]{Example}
\newtheorem{remark}[lemma]{Remark}

\numberwithin{equation}{section}

\begin{document}

\title{$\cD$-modules on Smooth Toric Varieties} 
 
\author[M. Musta\c{t}\v{a}]{Mircea~Musta\c{t}\v{a}}
\address{(Musta\c{t}\v{a}) Department of Mathematics \\ University of
California \\ Berkeley \\ CA \\ 94720--3840 \\ and Mathematics
Institute of the Romanian \newline Academy \\ {\tt
mustata@math.berkeley.edu}.}

\author[G. Smith]{Gregory~G.~Smith} 
\address{(Smith) Department of Mathematics \\ University of California
\\ \newline Berkeley \\ CA \\ 94720--3840 \\ {\tt
ggsmith@math.berkeley.edu}.}

\author[H. Tsai]{Harrison~Tsai}
\address{(Tsai) Department of Mathematics \\ University of California
\\ \newline Berkeley \\ CA \\ 94720--3840.}
\curraddr{Department of Mathematics \\ Cornell University \\ Ithaca \\
NY \\ 14853--4201 \\ {\tt htsai@math.cornell.edu}.}

\author[U. Walther]{Uli~Walther} 
\address{(Walther) Mathematical Sciences Research Institute \\
1000 Centennial Drive \\ Berkeley \\ CA \\ 94720--5070.}
\curraddr{Department of Mathematics \\ Purdue University \\ West
\newline Lafayette \\ IN \\ 47907--1395 \\ {\tt
walther@math.purdue.edu}.}

\date{15 July 2000}

\begin{abstract}
Let $X$ be a smooth toric variety.  Cox introduced the homogeneous
coordinate ring $S$ of $X$ and its irrelevant ideal $\mathfrak{b}$.
Let $A$ denote the ring of differential operators on $\Spec(S)$.  We
show that the category of $\cD$-modules on $X$ is equivalent to a
subcategory of graded $A$-modules modulo $\mathfrak{b}$-torsion.
Additionally, we prove that the characteristic variety of a
$\cD$-module is a geometric quotient of an open subset of the
characteristic variety of the associated $A$-module and that holonomic
$\cD$-modules correspond to holonomic $A$-modules.
\end{abstract}

\maketitle

\section{Introduction}

Let $X$ be a smooth toric variety over a field $k$.  Cox~\cite{cox}
introduced the homogeneous coordinate ring $S$ of $X$ and the
irrelevant ideal $\mathfrak{b}$.  The $k$-algebra $S$ is a polynomial
ring, with one variable for each one-dimensional cone in the fan
$\Delta$ defining $X$, and has a natural grading by the class group
$\Cl(X)$.  The monomial ideal $\mathfrak{b} \subset S$ encodes the
combinatorial structure of $\Delta$.  The following theorem of
Cox~\cite{cox} indicates the significance of the pair $(S,
\mathfrak{b})$.  We write $\text{{\sf $\cO$-Mod}}$ for the category of
quasi-coherent sheaves on $X$ and $\text{{\sf $S$-GrMod}}$ for the
category of graded $S$-modules.  A graded $S$-module $F$ is called
$\mathfrak{b}$-torsion if, for all $f \in F$, there exists $\ell > 0$
such that $\mathfrak{b}^{\ell}f = 0$.  Let $\btors$ denote the full
subcategory of $\mathfrak{b}$-torsion modules.

\begin{theorem*}[Cox] 
\begin{enumerate}
\item The category $\text{{\sf $\cO$-Mod}}$ is equivalent to the
quotient category $\text{{\sf $S$-GrMod}} / \btors$ .
\item The variety $X$ is a geometric quotient of $\Spec(S) \setminus
\Var(\mathfrak{b})$ by a suitable torus action.
\end{enumerate}
\end{theorem*}
\noindent When $X = \PP^{n}$, this is Serre's description of
quasi-coherent sheaves on projective space and the classical
construction of projective space.

The aim of this paper is to provide the $\cD$-module version of this
theorem --- $\cD$ denotes the sheaf of differential operators on $X$.
To state the analogue of the first part, we introduce the following
notation.  We write $\DMod$ for the category of left $\cD$-modules on
$X$.  The ring of differential operators on $\Spec(S)$ is the Weyl
algebra $A$; it also has a natural $\Cl(X)$ grading.  To each element
$\buu$ in $\Cl(X)^{\vee} = \Hom_{\ZZ}\big(\Cl(X), \ZZ \big)$, we
associate an ``Euler'' operator $\theta_{\buu} \in A$ [see
(\ref{eqn:Euler}) for the precise definition].  The full subcategory
of graded left $A$-modules $F$ satisfying $(\theta_{\buu} - \langle
\buu, \baa \rangle) \cdot F_{\baa} = 0$ for all $\baa \in \Cl(X)$ and
all $\buu \in \Cl(X)^{\vee}$ is denoted $\Ml$.

\begin{theorem}  \label{thm:mainresult} 
The quotient category $\Ml / \btors$ is equivalent to the category
$\DMod$.
\end{theorem}

\noindent The special case, when $X$ is projective space, can be found
in Borel~\cite{borel}.

This categorical equivalence is given by two functors.  The first
takes an object $F$ in $\Ml$ to the $\cD$-module $\tF$ whose sections
over the affine open subset $U_{\sigma}$ associated to $\sigma \in
\Delta$ are $(F_{x^{\hsigma}})_{\b00}$.  The second maps a
$\cD$-module $\cF$ to $ \Gamma_{L}(\cF) = \bigoplus_{\baa \in \Cl(X)}
H^{0}\big( X, \cO(\baa) \otimes \cF \big)$.  In fact, if $F$ is
finitely generated then $\tF$ is coherent and if $\cF$ is coherent
then it is of the form $\tF$ for some finitely generated graded
$A$-module $F$.  Our analysis of the $\Gamma_{L}(-)$ extends the work
of Musson~\cite{musson2} and Jones~\cite{jones2} on rings of twisted
differential operators on toric varieties.

Our second major result is

\begin{theorem} \label{thm:maincharvar}
If $F \in \Ml$ is finitely generated, then the characteristic variety
of $\tF$ is a geometric quotient of a suitable open subset of the
characteristic variety of $F$.
\end{theorem}

\noindent Moreover, given a finitely generated $F \in \Ml$ which has
no $\mathfrak{b}$-torsion, we show that the dimension of $\tF$ is
equal to the dimension of $F$ minus the rank of $\Cl(X)$.  In
particular, holonomic $A$-modules correspond to holonomic
$\cD$-modules.

The category of modules over the Weyl algebra is a well-studied
algebraic object and we hope to study $\cD$-modules on $X$ by using
these methods.  In particular, effective algorithms have been
developed for $\cD$-modules on affine space; for example, see the work
of Oaku~\cite{oaku}, Walther~\cite{walther},
Saito-Sturmfels-Takayama~\cite{sst} and Oaku-Takayama~\cite{ot}.  It
would interesting to use our results to extend these methods to smooth
toric varieties.

Many of our results are valid for a simplicial toric variety if one
replaces $S$ and $A$ with the subrings $\bigoplus_{\bbb \in \Pic(X)}
S_{\bbb}$ and $\bigoplus_{\bbb \in \Pic(X)} A_{\bbb}$; see
Cox~\cite{cox}.  This approach allows one to recover the general
results of Musson~\cite{musson2} and Jones~\cite{jones2}.  For
simplicity, we present the smooth case and leave the possible
generalizations to the reader.

The contents of this paper are as follows: The second section reviews
the basics about toric varieties, the Weyl algebra and $\cD$-modules.
In the third section, we determine the module associated with the
sheaf $\cD \otimes \cO(\bbb)$.  We introduce the $\Cl(X) \times
\Cl(X)$-graded $A$-$A$ bimodule
\[
\D = \bigoplus_{\bbb \in \Cl(X)} \frac{A(\bbb)}{A \cdot (\theta_{\buu}
+ \langle\buu, \bbb\rangle : \buu \in \Cl(X)^{\vee})}
\]
and construct a morphism
\[
\eta \colon \D \longrightarrow \bigoplus _{(\baa,\bbb)\in\Cl(X)^2}
H^0(X,\cO(\baa)\otimes\cD\otimes\cO(\bbb))
\]
We prove that $\eta$ is an isomorphism in two steps.  We first show
that $\eta$ induces an isomorphism on the associated sheaves.
Secondly, we establish that $H^{0}_{\mathfrak{b}}(\D) =
H^{1}_{\mathfrak{b}}(\D) = 0$.  The first step is a local statement
and can be reduced to results of Musson~\cite{musson1}.  We provide a
direct proof using simplifications due to Jones~\cite{jones1}.  The
fourth section contains the proof of Theorem~\ref{thm:mainresult}.  We
establish this result for both left and right $\cD$-modules.  In
general, there is an equivalence between these and we show how this is
induced at the level of $A$-modules.  In the last section, we prove
Theorem~\ref{thm:maincharvar} and related dimension results.

We would like to thank David Eisenbud, Bernd Sturmfels, Nobuki
Takayama and Will Traves for their helpful comments.

\section{Background}

We collect here a number of more or less standard definitions, results
and notation.  Throughout this paper, we work over an algebraically
closed field $k$ of characteristic zero.

\subsection*{Toric varieties}
Let $X$ be a smooth toric variety determined by the fan $\Delta$ in $N
\cong \ZZ^{n}$.  We write $\vv_{1}, \dotsc, \vv_{d}$ for the unique
lattice vectors generating the one-dimensional cones in $\Delta$ and
we assume that the $\vv_{i}$ span $N \otimes_{\ZZ} \RR$.  Each
$\vv_{i}$ corresponds to an irreducible torus invariant Weil divisor
in $X$.  Since these divisors generate the torus invariant Weil
divisors, we may identify the group of torus invariant Weil divisors
with $\ZZ^{d}$.  Let $\ee_{i}$ denote the standard basis for $\ZZ^{d}$
and set $\ee = \ee_{1} + \dotsb + \ee_{d}$.

There is a short exact sequence
\begin{equation} \label{eqn:ses}
0 \longrightarrow N^{\vee} \xrightarrow{\;\; \iota \;\;} \ZZ^{d}
\xrightarrow{\;\; \overline{(\;\;)} \;\;} \Cl(X) \longrightarrow 0 \, ,
\end{equation}
where $\iota(\pp) = \langle \pp, \vv_{1} \rangle \ee_{1} + \dotsb +
\langle \pp, \vv_{d} \rangle \ee_{d}$ and the second map $\aa \mapsto
\baa$ is the projection from Weil divisors to the divisor class group.
Since $X$ is smooth, the divisor class group $\Cl(X)$ is isomorphic to
the Picard group $\Pic(X)$.  In particular, the invertible sheaf (line
bundle) associated to $\baa \in \Cl(X)$ is denoted $\cO(\baa)$.

Following \cite{cox}, the homogeneous coordinate ring of $X$ is the
polynomial ring $S = k[x_{1}, \dotsc, x_{d}]$ with a $\Cl(X)$-grading
induced by
\[
\deg (x^{\aa}) = \deg (x_{1}^{\aa_{1}} \cdots x_{d}^{\aa_{d}}) = \baa \in
\Cl(X) \, .
\]
For a cone $\sigma \in \Delta$, $\hsigma$ is the set $\{ i : \vv_i
\not\in \sigma \}$ and $x^{\hsigma} = \prod_{\vv_{i} \not\in \sigma}
x_{i}$ is the associated monomial in $S$.  The irrelevant ideal of
$X$ is the reduced monomial ideal $\mathfrak{b} = (x^{\hsigma} :
\sigma\in\Delta)$.  

Recall that each cone $\sigma \in \Delta$ corresponds to an open
affine subset of $X$, $U_{\sigma} \cong \Spec
(S_{x^{\hsigma}})_{\b00}$.  Every graded $S$-module $F$ gives rise to
a quasi-coherent sheaf on $X$, denoted by $\tF$, which corresponds to
the module $(F_{x^{\hsigma}})_{\b00}$ over $U_{\sigma}$.  If $F$ is
finitely generated over $S$, then $\tF$ is a coherent $\cO$-module;
$\cO$ denotes the structure sheaf on $X$.  Moreover, every
quasi-coherent sheaf on $X$ is of the form $\tF$ for some graded
$S$-module $F$ and if the sheaf is coherent, then $F$ can be taken
finitely generated.  For an $S$-module $F$, we have $\tF = 0$ if and
only if $F = H^{0}_{\mathfrak{b}}(F)$; in other words $F$ is
$\mathfrak{b}$-torsion.

\subsection*{Weyl algebra}
By definition, the $d$-th Weyl algebra is 
\[
A = \frac{k \{ x_{1}, \dotsc, x_{d},\partial_{1}, \dotsc,
\partial_{d} \} }{\left( \begin{smallmatrix} x_{i}x_{j}-x_{j}x_{i}
& = & 0 \\ \partial_{i} \partial_{j} - \partial_{j} \partial_{i} & = &
0 \\ \partial_{i}x_{j} - x_{j}\partial_{i} & = & \delta_{ij}
\end{smallmatrix} \right)} \, .
\]
The canonical ring morphism $S \hookrightarrow A$ provides $A$ with
the structure of left $S$-module.  As in the case of $S$, $A$ has a
$\Cl(X)$-grading given by
\[
\deg (x^{\aa}\partial^{\bb}) = \deg (x_{1}^{\aa_{1}} \dotsm
x_{d}^{\aa_{d}}\partial_{1}^{\bb_{1}} \dotsm \partial_{d}^{\bb_{d}}) =
\baa - \bbb \in \Cl(X) \, ,
\]
and the $\baa$-th graded component of $A$ is denoted $A_{\baa}$.  For
each element $\buu$ of $\Cl(X)^{\vee} = \Hom_{\ZZ}\big( \Cl(X), \ZZ
\big)$, we have an Euler operator 
\begin{equation}
\label{eqn:Euler}
\theta_{\buu} = \langle \buu, \bee_{1} \rangle \theta_{1} + \dotsb +
\langle \buu, \bee_{d} \rangle \theta_{d} \, ,
\end{equation}
where $\theta_{i} = x_{i} \partial_{i}$.  Notice that $\theta_{\buu}$
has degree zero.

The Weyl algebra $A$ is isomorphic to the ring of differential
operators on $\AA^{d}$.  The natural action of $A$ on a polynomial $f
\in S$ is $x_{i} \bullet f = x_{i} \cdot f$ and $\partial_{i} \bullet
f = \frac{\partial f}{\partial x_{i}}$.  Since $S$ is also a subring
of $A$, the symbol $\bullet$ helps distinguish this action from the
product $\cdot \colon A \times A \to A$.  The ring of differential
operators on $S_{x^{\aa}} = S[x^{-\aa}]$ is denoted $A_{x^{\aa}}$ and
is equal to the localization $A[x^{-\aa}]$.

\subsection*{$\cD$-modules}

The sheaf of (algebraic) differential operators on $X$ is denoted
$\cD$.  On an affine open subset $U \subseteq X$, $H^{0}(U,\cD) =
\bigcup_{i \geq 0} \cD^{i}(U)$ where $\cD^{0}(U) = H^{0}(U,\cO)$ and
\[
\cD^{i}(U) : = \left\{ \text{$s \in \End_{k}\big( H^{0}(U, \cO) \big)$
: \begin{tabular}{l} $fs - sf \in \cD^{i-1}(U)$ \\ 
for all $f \in H^{0}(U,\cO)$ \end{tabular}}
\right\} \, .
\]
A $\cD$-module is a sheaf $\cF$ on $X$ which is quasi-coherent as an
$\cO$-module and has a structure of module over $\cD$. A $\cD$-module
is coherent if it is locally finitely generated over $\cD$.  We write
$\DMod$ and $\ModD$ for the categories of left and right $\cD$-modules
respectively.  The full subcategories of coherent left and right
$\cD$-modules are denoted $\DCoh$ and $\CohD$.

\section{Sheaves of Differential Operators} 

The goal of this section is to describe the left $S$-modules
corresponding to twists of the sheaf of differential operators.
Recall that, for a graded $A$-module $F$ and $\bbb \in \Cl(X)$,
$F(\bbb)$ is the shift of $F$ by $\bbb$: $F(\bbb)_{\baa} =
F_{\bbb+\baa}$.  We define the graded left $A$-modules
\[ 
\Dl(\bbb) = \frac{A(\bbb)}{A \cdot (\theta_{\buu} + \langle \buu,
\bbb \rangle : \buu \in \Cl(X)^{\vee})}
\]
and
\[
\D = \bigoplus_{\bbb \in \Cl(X)} \Dl(\bbb) \, .
\]
Notice that $\D$ has a $\Cl(X) \times \Cl(X)$-grading where
$\D_{(\baa,\bbb)} = \Dl(\bbb)_{\baa}$.

\begin{lemma}
The module $\D$ is a graded $A$-$A$ bimodule.
\end{lemma}

\begin{proof}
Multiplication in the ring $A$ yields the right action of $A$ on $\D$:
\[
\begin{CD}
A(\bbb)_{\baa} \otimes_{k} A_{\bbb'} @>>> A(\bbb + \bbb')_{\baa} \\
@VVV @VVV \\
\Dl(\bbb)_{\baa} \otimes_{k} A_{\bbb'} @>>> \Dl(\bbb + \bbb')_{\baa}
\, .
\end{CD}
\]
To see that the induced map is well-defined, observe that, for all
elements $f \in A_{\bbb'}$ and $\buu \in \Cl(X)^{\vee}$, we have
\begin{align*}
(\theta_{\buu} + \langle \buu, \bbb \rangle ) \cdot f 
&= f \cdot (\theta_{\buu} + \langle \buu, \bbb \rangle ) + \langle
\buu, \bbb' \rangle \cdot f \\
&= f \cdot (\theta_{\buu} + \langle \buu, \bbb + \bbb' \rangle ) \, .
\end{align*}
This action is clearly compatible with the left $A$-module structure.
It follows that $\D$ is an $A$-$A$ bimodule.  Since $f \in A_{\baa'}$
and $g \in A_{\bbb'}$ implies $f \cdot \D_{(\baa,\bbb)} \cdot g
\subseteq \D_{(\baa + \baa',\bbb + \bbb')}$, $\D$ is bigraded.  In other
words, if we let $A^{\circ}$ denote the opposite algebra, then $A
\otimes_{k} A^{\circ}$ is a $\Cl(X)^{2}$-graded ring and $\D$ is a
graded module over $A \otimes_{k} A^{\circ}$.
\end{proof}

Analogously, we define right $A$-modules:
\[ 
\Dr(\baa) = \frac{A(\baa)}{(\theta_{\buu} - \langle \buu, \baa \rangle
: \buu \in \Cl(X)^{\vee}) \cdot A}
\]
and
\[
\D' = \bigoplus_{\baa \in \Cl(X)} \Dr(\baa)\, .
\]
Again, $\D'$ is a $\Cl(X) \times \Cl(X)$-graded $A$-$A$ bimodule where
the multiplication on the left is induced by the multiplication in the
Weyl algebra.  In fact, we obtain the same module.

\begin{lemma} 
There is a canonical identification $\D = \D'$ which respects the graded
bimodule structure.
\end{lemma}

\begin{proof}
Since we have
\begin{align*}
\D_{(\baa, \bbb)} &= \frac{A_{\baa + \bbb}}{A_{\baa + \bbb} \cdot
(\theta_{\buu} + \langle \buu, \bbb \rangle : \buu \in \Cl(X)^{\vee})}
\intertext{and} 
\D_{(\baa, \bbb)}' &= \frac{A_{\baa + \bbb}}{ (\theta_{\buu} - \langle
\buu, \baa \rangle : \buu \in \Cl(X)^{\vee}) \cdot A_{\baa + \bbb}} \, ,
\end{align*}
it is enough to show that $A_{\baa + \bbb} \cdot (\theta_{\buu} +
\langle \buu, \bbb \rangle) = (\theta_{\buu} - \langle \buu, \baa
\rangle) \cdot A_{\baa+\bbb}$.  However, for every $f \in A_{\baa +
\bbb}$, we have
\begin{align*}
f \cdot (\theta_{\buu} + \langle \buu, \bbb \rangle) 
&= (\theta_{\buu} + \langle \buu, \bbb \rangle) \cdot f - \langle
\buu, \baa + \bbb \rangle \cdot f \\
&= (\theta_{\buu} - \langle \buu, \baa \rangle ) \cdot f \, ,
\end{align*}
which establishes the lemma.
\end{proof}

Now, the cohomology of $\cD$ has an $A$-$A$ bimodule structure.

\begin{lemma} 
The direct sum 
\[
\bigoplus_{(\baa,\bbb) \in \Cl(X)^{2}} H^{0}\big( X, \cO(\baa) \otimes
\cD \otimes \cO(\bbb) \big)
\]
is an $A$-$A$ bimodule.
\end{lemma}

\begin{proof}
It suffices to give $k$-linear maps:
\begin{multline} \label{eqn:mu_map}
\mu \colon A_{\baa'} \otimes_{k} H^{0}\big( X, \cO(\baa) \otimes \cD
\otimes \cO(\bbb) \big) \otimes_{k} A_{\bbb'} \\
\xrightarrow{\quad \quad} H^{0}\big( X, \cO(\baa+\baa') \otimes \cD
\otimes \cO(\bbb+\bbb') \big) \, .
\end{multline}
Locally, a section $s \in H^{0}\big(U_{\sigma},\cO(\baa) \otimes \cD
\otimes \cO(\bbb) \big)$ can be identified with an element of
$\Hom_{k}\big( (S_{x^{\hsigma}})_{-\bbb}, (S_{x^{\hsigma}})_{\baa})$,
where $\sigma$ is a cone in $\Delta$.  Moreover, the action of $f \in
A_{\baa}$ on $S$ descends to action on $S_{x^{\hsigma}}$ which
increases degrees by $\baa$.  Thus, we may define $(\mu
\vert_{U_{\sigma}})(f \otimes s \otimes g) = f \circ s \circ g$.  One
verifies that $\mu \vert_{U_{\sigma}}$ maps into $H^{0}\big(
U_{\sigma}, \cO(\baa+\baa') \otimes \cD \otimes \cO(\bbb+\bbb') \big)$
and that these local definitions glue to give the required map.
\end{proof}

We next construct a morphism of graded $A$-$A$ bimodules
\begin{equation} \label{eqn:eta_map}
\eta \colon \D \longrightarrow \bigoplus_{(\baa,\bbb)} H^{0}\big( X,
\cO(\baa) \otimes \cD \otimes \cO(\bbb) \big) \, .
\end{equation}
Since $\eta$ is a graded $k$-linear morphism, it is enough to define
$\eta(f)$ for $f \in \D_{(\baa, \bbb)}$.  For $f \in
\D_{(\baa,\bbb)}$, let $\eta(f)$ be the section of $H^{0}\big( X,
\cO(\baa) \otimes \cD \otimes \cO(\bbb) \big)$ whose restriction over
each $U_{\sigma}$ corresponds to the map induced by the action of $f
\colon (S_{x^{\hsigma}})_{-\bbb} \longrightarrow
(S_{x^{\hsigma}})_{\baa}$.  To see that $\eta(f)$ is well-defined,
consider $f \in A_{\baa + \bbb} \cdot (\theta_{\buu} + \langle \buu,
\bbb \rangle)$.  It follows that, for $g \in
(S_{x^{\hsigma}})_{-\bbb}$, we have
\[
(\theta_{\buu} + \langle \buu, \bbb \rangle ) \bullet g = - \langle
\buu, \bbb \rangle \cdot g + \langle \buu, \bbb \rangle \cdot g = 0
\]
and therefore $\eta(f) = 0$.  It is clear that $\eta$ is a morphism of
graded $A$-$A$ bimodules.  The main result of this section is the
following.

\begin{theorem} \label{thm:global_D}
The morphism $\eta$ $[$see equation~\eqref{eqn:eta_map}$]$ is an
isomorphism of graded $A$-$A$ bimodules.
\end{theorem}

Before proving Theorem~\ref{thm:global_D}, we collect some local results.
We first consider a local version of $\eta$.  By composing $\eta$ with
the restriction to $U_{\sigma}$ where $\sigma \in \Delta$, we obtain a
morphism of left $S$-modules
\begin{align*}
\D & \longrightarrow \bigoplus_{(\baa,\bbb) \in \Cl(X)^{2}} H^{0}\big(
U_{\sigma}, \cO(\baa) \otimes \cD \otimes \cO(\bbb) \big) \, ,
\intertext{which induces a morphism} 
\eta^{\sigma} \colon D_{x^{\hsigma}} & \longrightarrow
\bigoplus_{(\baa,\bbb) \in \Cl(X)^{2}} H^{0}\big( U_{\sigma},
\cO(\baa) \otimes \cD \otimes \cO(\bbb) \big) \, .
\end{align*}
Taking the degree zero component yields the ring morphism
\begin{equation} \label{eqn:varphi_map}
\bar{\varphi}_{\sigma} = \eta^{\sigma}_{(\b00, \b00)} \colon
\frac{(A_{x^{\hsigma}})_{\b00}}{(A_{x^{\hsigma}})_{\b00} \cdot
(\theta_{\buu} : \buu \in \Cl(X)^{\vee})} \longrightarrow
H^{0}(U_{\sigma},\cD) \, .
\end{equation}
Understanding $\bar{\varphi}_{\sigma}$ is an important ingredient in
establishing Theorem~\ref{thm:global_D}.  We begin by studying the map
$\varphi_{\sigma} \colon (A_{x^{\hsigma}})_{\b00} \longrightarrow
H^{0}(U_{\sigma}, \cD)$ which induces $\bar{\varphi}_{\sigma}$.

By definition, $H^{0}(U_{\sigma}, \cD)$ is the ring of differential
operators on the ring $(S_{x^{\hsigma}})_{\b00}$.  However, the
inclusion $\iota_{\sigma} \colon \sigma^{\vee} \cap N^{\vee}
\hookrightarrow \ZZ^{d}$ induces a ring isomorphism (denoted by the
same name) $\iota_{\sigma} \colon k[\sigma^{\vee} \cap N^{\vee}]
\xrightarrow{\cong} (S_{x^{\hsigma}})_{\b00}$.  To see this, observe
that $x^{\iota(\pp)} \in S_{x^{\hsigma}}$ if and only if $\langle \pp,
\vv_{i} \rangle \geq 0$, for all $\vv_{i} \in \sigma$, which is
equivalent to $\pp \in \sigma^{\vee}$.  It follows that the
isomorphism $\iota_{\sigma}$ induces an isomorphism, called
$\psi_{\sigma}$, from the differential operators on
$(S_{x^{\hsigma}})_{\b00}$ to the differential operators $R_{\sigma}$
on $k[\sigma^{\vee} \cap N^{\vee}]$.  More explicitly, we have
$\psi_{\sigma}(f) = \iota_{\sigma}^{-1} \circ f \circ \iota_{\sigma}$.
We will actually focus on the morphism $\psi_{\sigma} \circ
\varphi_{\sigma} \colon (A_{x^{\hsigma}})_{\b00} \longrightarrow
R_{\sigma}$.

Following \cite{musson1} and \cite{jones1}, we decompose
$(A_{x^{\hsigma}})_{\b00}$ and $R_{\sigma}$ under the appropriate
torus actions and express $\psi_{\sigma} \circ \varphi_{\sigma}$ in
terms of these decompositions.  To be more concrete, we identify
$N^{\vee}$ with $\ZZ^{n}$ by fixing a basis.  Let $\ep_{1}, \dotsc,
\ep_{n}$ denote the standard basis of $\ZZ^{n}$ and let $\ep = \ep_{1}
+ \dotsb + \ep_{n}$.  We continue to call the natural embedding $\iota
\colon N^{\vee} = \ZZ^{n} \longrightarrow \ZZ^{d}$.

The torus $(k^{\ast})^d$ acts on $k[x_{1}^{\pm 1}, \dotsc, x_{d}^{\pm
1}]$ by $\lambda * (x^\aa)=\lambda^\aa x^\aa$ which produces an action
on $A_{x^{\ee}}$ given by $\lambda * (x^{\aa}\partial^{\bb}) =
\lambda^{\aa - \bb}x^{\aa}\partial^{\bb}$ for $\lambda \in
(k^{\ast})^{d}$.  The corresponding eigenspace decomposition is $
A_{x^{\ee}} = \bigoplus_{\aa \in \ZZ^{d}} {x}^{\aa} \cdot W$, where
$W$ is the polynomial ring $k[\theta_{1}, \dotsc, \theta_{d}]$.
Taking the degree zero part, we have $(A_{x^{\ee}})_{\b00} =
\bigoplus_{\pp \in \ZZ^{n}} x^{\iota(\pp)} \cdot W$.  Since
$(A_{x^{\hsigma}})_{\b00}$ is invariant under the action of
$(k^{\ast})^{d}$, the decomposition of $(A_{x^{\ee}})_{\b00}$ yields
\[
(A_{x^{\hsigma}})_{\b00} = \bigoplus_{\pp \in \ZZ^{n}} x^{\iota(\pp)}
\cdot J(\pp) \, ,
\]
where $J(\pp)$ is an ideal in $W$.  To describe $J(\pp)$, recall that
$A_{x^{\hsigma}}$ is the ring of differential operators on
$S_{x^{\hsigma}}$.  Thus, if $g \in W$, then $x^{\iota(\pp)}g$ belongs
to $J(\pp)$ if and only if $(x^{\iota(\pp)}g) \bullet S_{x^{\hsigma}}
\subseteq S_{x^{\hsigma}}$.  Equivalently, for every $x^{\aa}$
satisfying $\aa_{i} \geq 0$ when $\vv_{i} \in \sigma$, we have
$g(\aa)x^{\iota(\pp)+\aa} \in S_{x^{\hsigma}}$.  We conclude that
$J(\pp)$ is the ideal of polynomials vanishing on
\[
Z(\pp) = \left\{ \aa \in \ZZ^{d} : \text{
\begin{minipage}[c]{150pt}
$\aa_{i} \geq 0$ if $\vv_{i} \in \sigma$ and \\
$\iota(\pp)_{j} + \aa_{j} < 0$ for some $\vv_j \in \sigma$
\end{minipage} 
} \right\} \, .
\]

Analogously, the affine space $\AA_{k}^{n}$ has as its associated Weyl
algebra,
\[
A' = \frac{k \{ y_{1}, \dotsc, y_{n}, \eth_{1}, \dotsc,
\eth_{n} \} }{\left( \begin{smallmatrix} y_{i}y_{j} - y_{j}y_{i}
& = & 0 \\ \eth_{i} \eth_{j} - \eth_{j} \eth_{i} & = &
0 \\ \eth_{i} y_{j} - y_{j}\eth_{i} & = & \delta_{ij}
\end{smallmatrix} \right)} 
\]
and the torus $(k^{\ast})^n$ acts on the Laurent ring $k[y_{1}^{\pm
1}, \ldots, y_{n}^{\pm 1}] = k[N^{\vee}]$ inducing an action on
$A_{y^{\ep}}'$.  The corresponding eigenspace decomposition is
$A_{y^{\ep}}' = \bigoplus_{\pp \in \ZZ^{n}} {y}^{\pp} \cdot W'$, where
$W'$ is the polynomial ring $k[\vartheta_{1}, \dotsc, \vartheta_{n}]$
and $\vartheta_{i} = y_{i}\eth_{i}$.  The inclusion $R_{\sigma}
\hookrightarrow A_{y^{\ep}}'$ identifies $R_{\sigma}$ with
\[
\left\{ f \in A_{y^{\ep}}' : f \bullet (k[\sigma^{\vee} \cap
N^{\vee}]) \subseteq k[\sigma^{\vee}\cap N^{\vee}] \right\} \, .
\] 
Since this condition is torus invariant, $R_{\sigma}$ is also
torus invariant and we obtain $R_{\sigma} = \bigoplus_{\pp \in
\ZZ^{n}} y^{\pp} \cdot I(\pp)$, where
\[
I(\pp) = \left\{ f \in W' : (y^{\pp} f) \bullet (k[\sigma^{\vee} \cap
N^{\vee}]) \subseteq k[\sigma^{\vee} \cap N^{\vee}] \right\} \, .
\]  
Identifying $W'$ with the coordinate ring of $\AA^{n}$, we have $f
\bullet y^{\qq} = f(\qq)y^{\qq}$ for every $\qq \in N^{\vee}$.  Hence,
$I(\pp)$ is the ideal of polynomials vanishing on
\[
Y(\pp) = \left\{ \qq \in \sigma^{\vee} \cap N^{\vee} : \qq + \pp
\not\in \sigma^{\vee} \cap N^{\vee} \right\} \, .
\]

Finally, we define a map $\rho \colon W \longrightarrow W'$.  The
inclusion $\iota \colon \ZZ^{n} \to \ZZ^{d}$ induces, by tensoring
with $k$, a linear embedding (denoted by the same name) $\iota :
\AA^{n} \to \AA^{d}$ of the corresponding affine spaces.  Identifying
$W$ and $W'$ with the coordinate rings of $\AA^{d}$ and $\AA^{n}$
respectively, we obtain the ring homomorphism $\rho = \iota^{\ast}
\colon W \to W'$.  Clearly, $\rho$ is surjective and $\Ker(\rho) =
\left( \theta_{\buu} : \buu \in \Cl(X)^{\vee} \right)$.

With this notation, we have

\begin{lemma} \label{lem:psi_phi}
The ring homomorphism $\psi_{\sigma} \circ \varphi_{\sigma} :
(A_{x^{\sigma}})_{\b00} \longrightarrow R_{\sigma}$ is given by
$x^{\iota(\pp)} \cdot g \mapsto y^{\pp} \cdot \rho(g)$ where $\pp \in
\ZZ^{n}$ and $g \in J(\pp)$.
\end{lemma}

\begin{proof}
It suffices to show that $(\psi_{\sigma} \circ \varphi_{\sigma})
(x^{\iota(\pp)} g)$ and $y^{\pp} \rho(g)$ have the same action on
$y^{\qq} \in k[\sigma^{\vee} \cap N^{\vee}]$.  On one hand, we have
\begin{align*}
(\psi_{\sigma} \circ \varphi_{\sigma})(x^{\iota(\pp)}g) \bullet
y^{\qq}
&= \iota^{-1} \left( \varphi_{\sigma}(x^{\iota(\pp)}g \bullet
\iota(y^{\qq}) \right) \\
&= \iota^{-1} \left( \varphi_{\sigma}(x^{\iota(\pp)}g \bullet
x^{\iota(\qq)} \right) \\
&= \iota^{-1} \left( g\big(\iota(\qq)\big) x^{\iota(\pp + \qq)}
\right) \\
&= g\big(\iota(\qq)\big) y^{\pp + \qq} \, .  
\end{align*}
On the other hand, we also have
\[
\big(y^{\pp}\rho(g)\big) \bullet y^{\qq} = \big( \rho(g)(\qq) y^{\pp +
\qq} \big) = g\big(\iota(\qq)\big)y^{\pp + \qq} \, ,
\]
which establishes the claim.
\end{proof}

Before returning our attention to $\bar{\varphi}_{\sigma}$, we need
one more lemma.

\begin{lemma} \label{lem:rho_theta}
If the elements $\theta_{i}$ and $\theta_{j}$ in $W = k[\theta_{1},
\dotsc, \theta_{d}]$ are distinct and correspond to rays $\vv_{i}$ and
$\vv_{j}$ in the same cone $\sigma$, then, for every pair of integers
$m_{i}$ and $m_{j}$, the elements $\rho(\theta_{i}) + m_{i}$ and
$\rho(\theta_{j}) + m_{j}$ are linearly independent over the field
$k$.
\end{lemma}

\begin{proof}
Suppose otherwise: for some $c \in k$, we have $\rho(\theta_{i}) +
m_{i} = c \cdot (\rho(\theta_{j}) + m_{j})$.  It follows that
$\rho(\theta_{i}) = c \cdot \rho(\theta_{j})$ and hence
\[
\sum_{\ell = 1}^{n} \langle \ep_{\ell}, \vv_{i} \rangle \vartheta_{\ell}
= c \cdot \left( \sum_{\ell = 1}^{n} \langle \ep_{\ell}, \vv_{j} \rangle
\vartheta_{\ell} \right) \, .
\]
We deduce that $\langle \ep_{\ell}, \vv_{i} \rangle = \langle
\ep_{\ell}, c \cdot \vv_{j} \rangle$ for all $\ell$ and thus $\vv_{i}
= c \cdot \vv_{j}$.  However, the $\vv_{i}$ in any cone $\sigma$ are
linearly independent.
\end{proof}

We are now in a position to understand $\bar{\varphi}_{\sigma}$.  In
particular, we obtain the following proposition which is a special
case of Musson's results on rings of differential operators; see
\cite{musson1}.

\begin{proposition}[Musson]  \label{prop:local_version}
For every $\sigma \in \Delta$, the map $\bar{\varphi}_{\sigma}$ $[$see
equation \eqref{eqn:varphi_map}$]$ is an isomorphism of rings.
\end{proposition}

\begin{proof}
The fact that $\bar{\varphi}_{\sigma}$ is a ring homomorphism follows
directly from the definition.  Thus, the assertion reduces to showing
that $\psi_{\sigma} \circ \varphi_{\sigma}$ is surjective and to
describe its kernel.  To achieve this, we determine the Zariski
closures of $Y(\pp)$ and $Z(\pp)$.  Since the Zariski closure of the
set of integer points inside a rational polyhedral cone is the linear
space spanned by that cone, it is easy to check that
\begin{equation} \label{eqn:Z&Y}
\begin{cases}
\overline{Z(\pp)} & = {\displaystyle \bigcup_{(i,m) \in \Lambda'}
\big\{ \bb \in k^{d} : b_i = m \big\},} \\
\overline{Y(\pp)} & = {\displaystyle \bigcup_{(i,m) \in \Lambda'}
\big\{ \qq \in k^{n} : \iota(\qq)_{i} = m \big\}} =
\iota^{-1}(\overline{Z(\pp)}) \, ,
\end{cases}
\end{equation}
where $\Lambda' = \{ (i,m) : \text{$\vv_{i} \in \sigma$ and $0 \leq m
\leq - \iota(\pp)_{i}$} \}$.  We claim that $\rho(J(\pp)) = I(\pp)$.
Indeed, equation \eqref{eqn:Z&Y} implies that $I(\pp) =
\sqrt{\rho(J(\pp))}$, so it is enough to check that $\rho(J(\pp))$ is
a radical ideal.  However, $J(\pp)$ is the principal ideal generated
by
\begin{equation} \label{eqn:h_p}
h_{\pp} : = \prod_{(i,m) \in \Lambda'} (\theta_{i} - m) \, .
\end{equation}
From the equation \eqref{eqn:h_p} and Lemma~\ref{lem:rho_theta}, we
see that $\rho(J(\pp))$ is reduced and hence $I(\pp) = \rho(J(\pp))$.
Applying Lemma~\ref{lem:psi_phi}, it follows that $\psi_{\sigma} \circ
\varphi_{\sigma}$ and $\bar{\varphi}_{\sigma}$ are surjective.  To
prove that $\psi_{\sigma} \circ \bar{\varphi}_{\sigma}$ is injective,
recall that $\Ker(\rho) = (\theta_{\buu} : \buu \in \Cl(X)^{\vee})$.
Thus, it is enough to show that, for every $\pp \in \ZZ^{n}$, we have
$J(\pp) \cap \Ker(\rho) = J(\pp) \cdot \Ker(\rho)$.  To see this,
observe that $f \in J(\pp) \cap \Ker(\rho)$ implies $f = h_{\pp}
f_{1}$.  Therefore, it suffices to notice that $\rho(h_{\pp}) \neq 0$,
which follows from Lemma~\ref{lem:rho_theta}.
\end{proof}

We now prove the main result in this section.

\begin{proof}[Proof of Theorem~\ref{thm:global_D}]
We must show that, for every $\bbb \in \Cl(X)$, the homomorphism of
graded left $A$-modules
\[
\eta_{(*,\bbb)} \colon \Dl(\bbb) \longrightarrow \bigoplus_{\baa \in
\Cl(X)} H^{0}\big( X, \cO(\baa) \otimes \cD \otimes \cO(\bbb) \big)
\] 
is an isomorphism.

The first step is to prove that $\eta_{(*,\bbb)}$ induces an
isomorphism of the associated sheaves.  Fix $\sigma \in \Delta$ and
choose $\bb \in \ZZ^{d}$, mapping to $\bbb$ in $\Cl(X)$, such that
$x^{\bb}$ is an invertible element in $S_{x^{\hsigma}}$.  Since the
restriction of $\cO(\bbb)$ to $U_{\sigma}$ is trivial, one may always
find such a $\bb$.  We then have a commutative diagram:
\[
\begin{CD}
\frac{(A_{x^{\hsigma}})_{\b00}}{(A_{x^{\hsigma}})_{\b00} \cdot
(\theta_{\buu}: \buu \in \Cl(X)^{\vee})} @>{\bar{\varphi}_{\sigma}}>>
H^0(U_{\sigma},\cD) \\
@V{\cdot x^{\bb}}VV @VV{ \otimes x^{\bb}}V \\
\frac{(A_{x^{\hsigma}})_{\bbb}}{(A_{x^{\hsigma}})_{\bbb}\cdot
(\theta_{\buu} + \langle \buu, \bbb \rangle :\buu \in \Cl(X)^{\vee})}
@>{\bar{\varphi}_{\sigma}'}>> H^{0}\big( U_{\sigma}, \cD \otimes
\cO(\bbb) \big)
\end{CD}
\]
where $\bar{\varphi}_{\sigma}$ is the morphism in
Proposition~\ref{prop:local_version} and $\bar{\varphi}_{\sigma}'$ is
the analogous morphism induced by $\eta_{(\b00,\bbb)}$.  Now,
Proposition~\ref{prop:local_version} implies that
$\bar{\varphi}_{\sigma}$ is an isomorphism and the vertical arrows are
clearly isomorphisms.  It follows that $\bar{\varphi}_{\sigma}'$ is an
isomorphism and therefore $\eta_{(*,\bbb)}$ induces an isomorphism of
the associated sheaves.

If $F$ is a graded $S$-module, we write
\[ 
\Gamma_{L}(\tF) = \bigoplus_{\baa \in \Cl(X)} H^{0}\big(X, \cO(\baa)
\otimes \tF \big) \, .
\]
For every such $F$, there is an exact sequence
\begin{equation} \label{eqn:local_global}
0 \longrightarrow H^{0}_{\mathfrak{b}}(F)
\longrightarrow F \longrightarrow \Gamma_{L}(\tF) \longrightarrow
H^{1}_{\mathfrak{b}}(F) \longrightarrow 0 \, ,
\end{equation}
where $\mathfrak{b}$ is the irrelevant ideal; see \cite{ems}.  Hence,
if $H^{0}_{\mathfrak{b}} \big( \Dl(\bbb) \big)$ and
$H^{1}_{\mathfrak{b}} \big( \Dl(\bbb) \big)$ both vanish, then
$\eta_{(*,\bbb)}$ is an isomorphism.  We relegated these vanishing
results to Propositions~\ref{prop:nonzerodivisor} and \ref{prop:H1=0}
below.
\end{proof}

Our first vanishing result is

\begin{proposition} \label{prop:nonzerodivisor}
The element $x^{\ee} \in \mathfrak{b}$ is not a zero divisor on
$\Dl(\bbb)$ and $H^{0}_{\mathfrak{b}}(\Dl(\bbb)) = 0$.
\end{proposition}

\begin{proof}
The first assertion implies the second, so it suffices to show that
$x^{\ee}$ is not a zero divisor.  Every $\aa \in \ZZ^{d}$ can be
written uniquely as $\aa = \aa^{+} - \aa^{-}$ where $\aa^{+}$ and
$\aa^{-}$ are non-negative and have disjoint support.  Consider the
action of the torus $(k^{*})^d$ on the Weyl algebra $A$; the
corresponding eigenspace decomposition is $A = \bigoplus_{\aa \in
\ZZ^{d}} x^{\aa^{+}}\partial^{\aa^{-}} \cdot W$, where $W =
k[\theta_{1}, \dotsc, \theta_{d}]$.  Let $L_{0} \subset W$ be the
ideal generated by $\theta_{\buu} + \langle \buu, \bbb \rangle$ for
all $\buu \in \Cl(X)^{\vee}$.  Since $L_{0}$ is generated by linear
forms, it is a prime ideal.  Let $L$ denote the left $A$-ideal
$(\theta_{\buu} + \langle \buu, \bbb \rangle : \buu \in
\Cl(X)^{\vee})$.  With this notation, we have $L = \bigoplus_{\aa \in
\ZZ^{d}} x^{\aa^{+}}\partial^{\aa^{-}} \cdot L_{0}$.

For $x^{\ee}$ to be a non-zero divisor on $\Dl(\bbb)$, it suffices to
prove that, for $g \in W$ and $\aa \in \ZZ^{d}$, the relation $x^{\ee}
x^{\aa^{+}}\partial^{\aa^{-}} g \in L$ implies $g \in L_{0}$.  To
accomplish this, we first note
\begin{align*}
x^{\ee} \cdot x^{\aa^{+}}\partial^{\aa^{-}} &= \left( \prod_{\aa_{i}
\geq 0} x_{i}^{\aa_{i}+1} \right) \left( \prod_{\aa_{i} < 0} x_{i}
\partial^{-\aa_{i}} \right) \\ 
&= \left( \prod_{\aa_{i} \geq 0} x_{i}^{\aa_{i}+1} \right) \left(
\prod_{\aa_{i} < 0} \partial^{-\aa_{i}-1} (\theta_{i} + \aa_{i} + 1)
\right) \\
&= x^{(\aa + \ee)^{+}}\partial^{(\aa + \ee)^{-}} \prod_{\aa_{i} < 0}
(\theta_{i} + \aa_{i} + 1) \, .
\end{align*}
Hence, $x^\ee \cdot x^{\aa^{+}}\partial^{\aa^{-}} g \in L$ implies
$\left( \prod_{\aa_{i} < 0} (\theta_{i} + \aa_{i} + 1) \right) \cdot g
\in L_{0}$.  Suppose $g \not\in L_{0}$.  Since $L_{0}$ is a prime
ideal, there exists an index $i$ such that $\aa_{i} < 0$ and
$\theta_{i} + \aa_{i} + 1 \in L_{0}$.  Expressing $\theta_{i} +
\aa_{i} + 1$ in terms of the linear generators of $L_{0}$, it follows
that there is $\buu \in \Cl(X)^{\vee}$ such that $\theta_{i} =
\theta_{\buu}$.  However, for every $\ww \in N^{\vee}$,
equation~\eqref{eqn:ses} implies
\[
\big\langle \buu, \langle \ww, \vv_{1} \rangle \bee_{1} + \dotsb +
\langle \ww, \vv_{d} \rangle \bee_{d} \big\rangle = 0 \, ,
\]  
from which we deduce that $\langle \ww, \vv_{i} \rangle = 0$ and
$\vv_{i} = 0$ giving a contradiction.
\end{proof}

We end this section with

\begin{proposition} \label{prop:H1=0}
The local cohomology module $H^{1}_{\mathfrak{b}}(\Dl(\bbb))$
vanishes.
\end{proposition}

\begin{proof}
Since $x^{\ee}$ is not a zero divisor on $\Dl(\bbb)$, there is a short
exact sequence of $S$-modules
\[
0 \longrightarrow \Dl(\bbb) \xrightarrow{\;\; x^{\ee} \cdot \;\;}
\Dl(\bbb) \longrightarrow Q = \frac{\Dl(\bbb)}{x^{\ee} \cdot
\Dl(\bbb)} \longrightarrow 0 \, ,
\]
and the long exact sequence of local cohomology gives 
\begin{equation} \label{eqn:H1}
0 \longrightarrow H^{0}_{\mathfrak{b}}(Q) \longrightarrow
H^{1}_{\mathfrak{b}}(\Dl(\bbb)) \xrightarrow{\;\; x^{\ee} \cdot \;\;}
H^{1}_{\mathfrak{b}}(\Dl(\bbb)) \, .
\end{equation}
Because every element in $H^{1}_{\mathfrak{b}}\big( \Dl(\bbb) \big)$
is annihilated by a power of $\mathfrak{b}$, the injectivity of
$x^{\ee} \cdot$ in equation~\eqref{eqn:H1} implies that
$H^{1}_{\mathfrak{b}} \big( \Dl(\bbb) \big) = 0$.  Thus, it suffices
to prove that $ H^{0}_{\mathfrak{b}}(Q) = 0$

Let $K$ be the left $A$-ideal satisfying $Q = A(\bbb)/K$.  To prove
that $H^{0}_{\mathfrak{b}}(Q) = 0$, we must show that if $f \in A$
satisfies $(x^{\hsigma})^m f \in K$ for some $m \geq 1$ and all
$\sigma \in \Delta$, then $f \in K$.  Using the notation from the
proof of Proposition~\ref{prop:nonzerodivisor}, we have $K = L +
x^{\ee} \cdot A$.  From the decomposition of $A$, we have
\begin{align*}
x^{\ee} \cdot A &= \bigoplus_{\aa \in \ZZ^{d}} x^\ee \cdot x^{\aa^{+}}
\partial^{\aa^{-}} \cdot W \\ 
&= \bigoplus_{\aa \in \ZZ^{d}} \left( x^{(\aa+\ee)^{+}}
\partial^{(\aa+\ee)^{-}} \left(\prod_{\aa_{i} < 0} (\theta_{i} +
\aa_{i} + 1) \right) \cdot W \right) \, ,
\end{align*}
and we deduce $K = \bigoplus_{\aa \in \ZZ^{d}}
x^{\aa^{+}}\partial^{\aa^{-}} \cdot K(\aa)$, where $K(\aa)$ is defined
to be $L_{0} + \big( \prod_{\aa_{i} \leq 0} (\theta_{i} + \aa_{i})
\big) \cdot W$.  Thus, it is enough to consider elements $f \in A$ of
the form $x^{\aa^{+}}\partial^{\aa^{-}}g$ with $g \in W$ and prove
that $g \in K(\aa)$.  Moreover, we may assume that $m + \aa_{i} > 0$
for $1 \leq i \leq d$.

By induction on $r$, we see that $x_{i}^{r}\partial_{i}^{r} =
\prod_{j=1}^{r}(\theta_{i} -j+1)$ for $r \geq 1$.  Hence, for $\aa_{i}
< 0$, we have
\[
x_{i}^{m} \partial^{-\aa_{i}} = x_{i}^{m + \aa_i}x_{i}^{-\aa_{i}}
\partial_{i}^{-\aa_{i}} = x_{i}^{m+\aa_{i}} \cdot
\prod_{j=1}^{-\aa_{i}}(\theta_{i}-j+1)
\]
and, for $\sigma \in \Delta$, we obtain
\[
(x^{\hsigma})^{m} x^{\aa^{+}}\partial^{\aa^{-}} = \left( x^{(\aa +
m\hsigma)^{+}}\partial^{(\aa + m \hsigma)^{-}} \right) \cdot \prod_{i
\in \widehat{\Lambda}} \prod_{j=1}^{-\aa_{i}}(\theta_{i}-j+1) \, ,
\]
where $\widehat{\Lambda} = \{ i : \text{$\vv_{i} \notin \sigma$ and
$\aa_{i} < 0$} \}$.  We deduce
\begin{equation} \label{eqn:K=L+}
\left( \prod_{i \in \widehat{\Lambda}}
\prod_{j=1}^{-\aa_{i}}(\theta_{i}-j+1) \right) \cdot g \in L_{0} +
\left( \prod_{i \in \Lambda} (\theta_{i} + \aa_{i}) \right) \cdot W \,
,
\end{equation}
where ${\Lambda} = \{ i : \text{$\vv_{i} \in \sigma$ and $\aa_{i} \leq
0$} \}$.

For each $\bb \in \ZZ^{d}$, we define an automorphism $\alpha_{\bb}
\colon W \longrightarrow W$ given by $\alpha_{\bb}(\theta_{i}) =
\theta_{i} - \bb_{i}$ and we define $\rho_{\bb} \colon W
\longrightarrow W'$ to be the composition $\rho_{\bb} = \rho \circ
\alpha_{\bb}$ (the map $\rho$ is defined in the paragraph before
Lemma~\ref{lem:psi_phi}).  It is clear that $\rho_{\bb}$ is surjective
and $\Ker(\rho_{\bb}) = L_{0}$.  Applying $\rho_{\bb}$ to
equation~\eqref{eqn:K=L+} gives
\[
\left( \prod_{i \in \widehat{\Lambda}} \prod_{j=1}^{-\aa_{i}}
\rho_{\bb}(\theta_{i}-j+1) \right) \cdot \rho_{\bb}(g) \in \left(
\prod_{i \in {\Lambda}} \rho_{\bb}(\theta_{i} + \aa_{i}) \right) \cdot
W' \, .
\]
Lemma~\ref{lem:rho_theta} obviously extends to $\rho_{\bb}$ and
implies
\[
\rho_{\bb}(g) \in \left( \prod_{i \in {\Lambda}} \rho_{\bb}(\theta_{i}
+ \aa_{i}) \right) \cdot W' \, .
\]
Since this relation holds for every $\sigma \in \Delta$, a second
application of Lemma~\ref{lem:rho_theta} shows that
$\rho_{\bb}(g) \in \big( \prod_{\aa_{i} \leq 0} \rho_{\bb}(\theta_{i}
+ \aa_{i}) \big) \cdot W'$ and therefore $g \in K(\aa) = L_{0} + \big(
\prod_{a_{i} \leq 0} (\theta_{i} + \aa_{i}) \big) \cdot W$.
\end{proof}

\section{$\cD$-modules} \label{section:D_modules}

We now use Theorem~\ref{prop:local_version} to describe the relation
between $A$-modules and $\cD$-modules on $X$.  We begin by showing
that the $S$-module associated to a $\cD$-module has a graded
$A$-module structure.

\begin{proposition} \label{prop:grmod}
If $\cF$ is a left $\cD$-module, then the graded $S$-module
\begin{align*}
\Gamma_{L}(\cF) &= \bigoplus_{\baa \in \Cl(X)} H^{0}\big( X, \cO(\baa)
\otimes \cF \big) \, ,
\intertext{has a graded left $A$-module structure, extending the left
$S$-module structure.  Similarly, if $\cG$ is a right $\cD$-module,
then the graded $S$-module}
\Gamma_{R}(\cG) &= \bigoplus_{\bbb \in \Cl(X)} H^{0}\big( X, \cG
\otimes \cO(\bbb)\big) \, ,
\end{align*}
has a graded right $A$-module structure, extending the right
$S$-module structure.
\end{proposition}

\begin{proof}
We present the left $\cD$-modules case here --- the proof for right
$\cD$-modules is completely analogous.  For the first assertion, it is
enough as well to construct $k$-linear maps
\[
\mu_{\baa',\baa}^{\cF} \colon A_{\baa'} \otimes_{k} H^{0}\big( X,
\cO(\baa) \otimes \cF \big) \longrightarrow H^{0}\big(X,
\cO(\baa+\baa') \otimes \cF \big) \, ,
\] 
for all $\baa, \baa' \in \Cl(X)$, satisfying the obvious axioms.  To
accomplish this, we consider a local version of the left
multiplication map \eqref{eqn:mu_map} when $\bbb = 0$.  More
explicitly, for each $\sigma \in \Delta$, this morphism is
\[
\mu_{\baa', \baa} \vert_{U_{\sigma}} \colon A_{\baa'} \otimes_{k}
H^{0}\big(U_{\sigma}, \cO(\baa) \otimes \cD \big) \longrightarrow
H^{0}\big(U_{\sigma}, \cO(\baa+\baa') \otimes \cD \big) \, ,
\]
given by $\mu_{\baa',\baa} \vert_{U_{\sigma}}(f \otimes s) = f \circ
s$.  We claim that $\mu_{\baa', \baa} \vert_{U_{\sigma}}$ is a
morphism of right $H^{0}(U_{\sigma}, \cD)$-modules.  To see this,
recall that sections $s \in H^{0}(U_{\sigma}, \cO(\bbb) \otimes \cD)$
and $g \in H^{0}(U_{\sigma}, \cD)$ can be identified with elements in
$\Hom_{k}\big( (S_{x^{\hsigma}})_{0}, (S_{x^{\hsigma}})_{\bbb} \big)$
and $\Hom_{k}\big( (S_{x^{\hsigma}})_{0}, (S_{x^{\hsigma}})_{0} \big)$
respectively.  In particular, we have
\[
\big(\mu_{\baa', \baa} \vert_{U_{\sigma}}(f \otimes s)\big) \cdot g =
f \circ s \circ g = \mu_{\baa', \baa} \vert_{U_{\sigma}}(f \otimes s
\cdot g) \, ,
\]
for all $f \in A_{\baa'}$.  It follows that, by tensoring the map
$\mu_{\baa', \baa} \vert_{U_{\sigma}}$ on the right with
$H^{0}(U_{\sigma}, \cF)$ over $H^{0}(U_{\sigma}, \cD)$, we obtain a
$k$-linear map
\[
\mu_{\baa',\baa}^{\cF}|_{U_{\sigma}} \colon A_{\baa'} \otimes_{k}
H^{0}\big( U_{\sigma}, \cO(\baa) \otimes \cF \big) \longrightarrow
H^{0}\big(U_{\sigma}, \cO(\baa'+\baa) \otimes \cF \big) \, .
\]
These maps glue together to give $\mu_{\baa',\baa}^{\cF}$ which makes
$F$ into a graded left $A$-module.
\end{proof}

Let $\Ml$ be the category of graded left $A$-modules $F$ such that
\begin{equation} \label{eqn:theta_condition}
(\theta_{\buu} - \langle \buu, \baa \rangle) \cdot F_{\baa} = 0 \text{
for all $\baa \in \Cl(X)$ and all $\buu \in \Cl(X)^{\vee}$.}
\end{equation}
A graded $A$-module $F$ is called $\mathfrak{b}$-torsion if, for every
$f \in F$, there exists $\ell > 0$ such that $\mathfrak{b}^{\ell}f =
0$.  Let $\btors$ denote the full subcategory of
$\mathfrak{b}$-torsion modules.  Similarly, $\Mr$ is the category of
graded right $A$-modules $G$ satisfying $G_{\bbb} \cdot (\theta_{\buu}
+ \langle \buu, \bbb \rangle) = 0$ for all $\bbb \in \Cl(X)$ and all
$\buu \in \Cl(X)^{\vee}$.  Let $\torsb$ denote the full subcategory of
$\mathfrak{b}$-torsion modules in $\Mr$.  It is clear that $\Ml$ and
$\Mr$ are both abelian categories closed under taking graded
subquotients.  The main result in this section is:

\begin{theorem} \label{thm:Dmodules}
The map $F \mapsto \tF$ is an exact functor from $\Ml$ to $\DMod$,
$\cF \mapsto \Gamma_{L}(\cF)$ is a left exact functor from $\DMod$ to
$\Ml$ and there are natural transformations of functors 
\[
{\rm id}_{\Ml / \btors } \xrightarrow{\;\; \cong \;\;} \Gamma_{L} \circ
\widetilde{\quad} 
\quad \text{ and } \quad 
\widetilde{\quad} \circ \Gamma_{L} \xrightarrow{\;\; \cong \;\;} {\rm
id}_{\DMod} \, .
\]

Similarly, the map $G \mapsto \tG$ is an exact functor from $\Mr$ to
$\ModD$, $\cG \mapsto \Gamma_{R}(\cG)$ is a left exact functor from
$\ModD$ to $\Mr$ and there are natural transformations of functors
\[
{\rm id}_{\Mr / \torsb} \xrightarrow{\;\; \cong \;\;} \Gamma_{R} \circ
\widetilde{\quad}
\quad \text{ and } \quad 
\widetilde{\quad} \circ \Gamma_{R} \xrightarrow{\;\; \cong \;\;} {\rm
id}_{\ModD} \, .
\]
\end{theorem}

In particular, every left $\cD$-module is of the form $\tF$ for some
graded left $A$-module $F$ and every right $\cD$-module is of the form
$\tG$ for some graded right $A$-module $G$.

\begin{proof}[Proof of Theorem~\ref{thm:Dmodules}]
Again, we give the proof only for left modules.  For the first part,
we consider an object $F$ in $\Ml$.  By definition, we have
$H^{0}(U_{\sigma}, \tF) = (F_{x^{\hsigma}})_{\b00}$, where
$(F_{x^{\hsigma}})_{\b00}$ is a left $(A_{x^{\hsigma}})_{\b00}$-module
and $\sigma \in \Delta$.  In light of
Theorem~\ref{prop:local_version}, we must show that, for every $\buu
\in \Cl(X)^{\vee}$, we have $\theta_{\buu} \cdot
(F_{x^{\hsigma}})_{\b00} = 0$.  Now, if $\frac{f}{(x^{\hsigma})^{m}}
\in (F_{x^{\hsigma}})_{\b00}$, then $f \in F_{m\bee_{\sigma}}$ where
$\ee_{\sigma} = \sum_{i \in \hsigma} \ee_{i}$ and the hypotheses on
$F$ imply
\begin{align*}
\theta_{\buu} \cdot \tfrac{f}{(x^{\hsigma})^m} &= \left( \theta_{\buu}
\tfrac{1}{(x^{\hsigma})^m} \right) \cdot f \\ 
&= \tfrac{1}{(x^{\hsigma})^m} \theta_{\buu} \cdot f - \sum_{i \in
\hsigma} m \langle \buu, \bee_{i} \rangle \tfrac{1}{(x^{\hsigma})^m}
\cdot f \\
&= \tfrac{1}{(x^{\hsigma})^m} \big(\theta_{\buu} - \langle \buu,
m\bee_{\sigma} \rangle \big) \cdot f = 0 \, .
\end{align*}
Therefore $H^{0}(U_{\sigma}, \tF)$ has a structure of left module over
$H^{0}(U_{\sigma}, \cD)$ for every $\sigma \in \Delta$.  It is
straightforward to verify that these structures glue together to give
a $\cD$-module structure on $\tF$.  By construction, we see that $\tF$
is quasi-coherent sheaf over $\cD$.

Conversely, let $\cF$ be an object of $\DMod$.  Applying
Proposition~\ref{prop:grmod}, we know that $F = \Gamma_{L}(\cF)$ is a
graded left $A$-module, so it is enough to prove that $F$ satisfies
\eqref{eqn:theta_condition}.  Fixing $\buu \in \Cl(X)^{\vee}$ and
$\baa \in \Cl(X)$, it suffices to show that
\[
\mu_{\b00, \baa}^{\cF}|_{U_{\sigma}}\big( (\theta_{\buu} - \langle
\buu, \baa \rangle ) \otimes s' \big) = 0 \, ,
\] 
for every $\sigma \in \Delta$ and all sections $s' \in H^{0}\big(
U_{\sigma}, \cO(\baa) \otimes \cF \big)$.  As explained in
Proposition~\ref{prop:grmod}, we have
\begin{align*}
H^{0} \big( U_{\sigma}, \cO(\baa) \otimes \cF \big) &= H^{0} \big(
U_{\sigma}, \cO(\baa) \otimes \cD \otimes_{\cD} \cF \big) \\
&= H^{0} \big( U_{\sigma}, \cO(\baa) \otimes \cD \big)
\otimes_{H^{0}(U_{\sigma}, \cD)} H^{0}(U_{\sigma}, \cF) \, ,
\end{align*}
so that we may identify $s'$ with a linear combination of elements of
the form $s \otimes f$.  By definition, we have
\[
\mu_{\b00, \baa}^{\cF}|_{U_{\sigma}} \big( (\theta_{\buu} - \langle
\buu, \baa \rangle) \otimes s \otimes f \big) = (\theta_{\buu} -
\langle \buu, \baa \rangle) \circ s \otimes f\, ,
\]
and we claim that $(\theta_{\buu} - \langle \buu, \baa \rangle) \circ
s$ is zero.  Indeed, for every $x^{\cc} \in (S_{x^{\hsigma}})_{\baa}$,
we have $ \theta_{\buu}(x^{\cc}) = \langle \buu, \bee_{i} \rangle
c_{1} x^{\cc} + \dotsb + \langle \buu, \bee_{i} \rangle c_{d} x^{\cc}
= \langle \buu, \baa \rangle x^{\cc}$.  Therefore, if $g \in
H^{0}(U_{\sigma},\cO) = (S_{x^{\hsigma}})_{0}$, then $s(g) \in
(S_{x^{\hsigma}})_{\baa}$ and
\[
(\theta_{\buu} - \langle \buu, \baa \rangle)(s(g)) = \langle \buu,
\baa \rangle s(g) - \langle \buu, \baa \rangle s(g) = 0 \, .
\]

Finally, the exact sequence \eqref{eqn:local_global} provides the
first natural transformation, once we observe that
$H_{\mathfrak{b}}^{i}(F)$ is $\mathfrak{b}$-torsion.  It follows from
Cox~\cite{cox} that the sheaf associated to $\Gamma_{L}(\cF)$ is
isomorphic to $\cF$.
\end{proof}

As a corollary, we obtain 

\begin{proof}[Proof of Theorem~\ref{thm:mainresult}]
Follows immediately from Theorem~\ref{thm:Dmodules}.
\end{proof}

We next turn our attention to coherent $\cD$-modules and finitely
generated graded $A$-modules.  We write $\Mlf$ and $\Mrf$ for the full
subcategories of $\Ml$ and $\Mr$ consisting of finitely generated
$A$-modules.

\begin{proposition}
If $F$ is an object in $\Mlf$ then $\tF$ is a coherent left
$\cD$-module.  Moreover, every coherent left $\cD$-module is of the
form $\tF$ for some $F \in \Mlf$.  Similarly, $G \in \Mrf$ implies
$\tG$ belongs to $\CohD$ and every coherent right $\cD$-module is
isomorphic to $\tG$ for some $G \in \Mrf$.
\end{proposition}

Our proof is analogous to the $\cO$-module case found in
Cox~\cite{cox}.

\begin{proof}
Once again, we present the proof only for left modules.  Suppose that
$F$ belongs to $\Mlf$.  To establish that $\tF \in \DCoh$, we have to
check that, for every $\sigma \in \Delta$, $(F_{x^{\hsigma}})_{\b00}$
is finitely generated over $(A_{x^{\hsigma}})_{\b00}$.  Thus, it
suffices to show that, for every element $f \in F$, there is an
invertible element $g \in S_{x^{\hsigma}}$ such that $g \cdot f$ has
degree zero.  Consider $f$ in $F_{\baa}$.  Since $X$ is smooth, there
is a divisor corresponding to $\aa \in \ZZ^{d}$ supported outside
$\sigma$ and whose class is $\baa$.  It follows that $x^{-\aa}$ is an
invertible element in $S_{x^{\hsigma}}$ and $x^{-\aa} \cdot f$ has
degree zero.

For the second assertion, we must show that given a coherent left
$\cD$-module $\cF$ there exists a finitely generated $A$-submodule $F$
of $\Gamma_{L}(\cF)$ such that $\widetilde{F} \cong \cF$.  Since $\cF$
is coherent, $H^{0}(U_{\sigma}, \cF)$ is finitely generated over
$H^{0}(U_{\sigma},\cD)$, for every $\sigma \in \Delta$.  Choose, for
each $\sigma \in \Delta$, a finite set of homogeneous elements in
$\Gamma_{L}(\cF)$ which are the numerators for a corresponding set of
generators of $H^{0}(U_{\sigma}, \cF)$.  Setting $F$ to be the
$A$-submodule of $\Gamma_{L}(\cF)$ generated by the union of these
sets, we have $\tF \cong \cF$ and $F$ is finitely generated over $A$.
\end{proof}

\begin{remark} \label{rem:rep}
As a consequence of Theorem~\ref{thm:Dmodules}, we have $\tF = 0$ if
and only if the graded $S$-module $F$ satisfies $F =
H_{\mathfrak{b}}^{0}(F)$; this is equivalent to saying that $F$ is a
$\mathfrak{b}$-torsion module.  At the other extreme, $F$ has no
$\mathfrak{b}$-torsion when $H^{0}_{\mathfrak{b}}(F) = 0$ and we say
that $F$ is $\mathfrak{b}$-saturated when $H^{0}_{\mathfrak{b}}(F) =
H^{1}_{\mathfrak{b}}(F) = 0$.  Now, every left $\cD$-module $\cF$ can
be represented by a unique saturated $A$-module, namely
$\Gamma_{L}(\cF)$.  Unfortunately, this may not be finitely generated,
even if $\cF$ is coherent.  However, by replacing $F$ with a suitable
submodule of $F / H^{0}_{\mathfrak{b}}(F)$, we may assume that $F$ has
no $\mathfrak{b}$-torsion and is finitely generated, whenever $\cF$ is
coherent.
\end{remark}

\begin{example}
The $A$-module corresponding to the structure sheaf $\cO$ (which is a
left $\cD$-module) is $\Gamma_{L}(\cO) = \frac{A}{A \cdot
(\partial_{1}, \dotsc, \partial_{d})}$.  In particular,
$\Gamma_{L}(\cO)$ is isomorphic to $S$, where $S$ has the standard
$A$-module structure.
\end{example}

\begin{corollary} \label{cor:twisted_D}
For every $\bbb \in \Cl(X)$, there is an isomorphism of graded left
$A$-modules $\Gamma_{L}\big( \cD \otimes \cO(\bbb) \big) \cong
\Dl(\bbb)$.  Similarly, we have $\Gamma_{R} \big( \cO(\baa) \otimes
\cD \big) \cong \Dr(\baa)$.
\end{corollary}

\begin{proof}
Theorem~\ref{thm:global_D} provides $\eta_{(*,\b00)} \colon \Dl(\bbb)
\xrightarrow{\;\; \cong \;\;} \Gamma_{L}\big( \cD \otimes \cO(\bbb)
\big)$ and $\eta_{(\b00,*)} \colon \Dr(\baa) \xrightarrow{\;\; \cong
\;\;} \Gamma_{R}\big( \cO(\baa) \otimes \cD \big)$.
\end{proof}

\begin{corollary} \label{cor:D_in_M}
For $\baa, \bbb \in \Cl(X)$, we have $\Dl(\baa) \in \Ml$ and
$\Dr(\bbb) \in \Mr$.
\end{corollary}

\begin{proof}
This follows from Corollary~\ref{cor:twisted_D} and
Theorem~\ref{thm:Dmodules}.
\end{proof}

\begin{corollary} \label{cor:presentation1}
Let $F$ be a graded left $A$-module generated by homogeneous elements
$\{f_{i} \}_{i}$ with $\deg f_{i} = \bbb_{i}$.  For $F$ to belong to
$\Ml$, it is necessary and sufficient that $(\theta_{\buu} - \langle
\buu ,\bbb_{i} \rangle) \cdot f_{i} = 0$ for all $i$ and all $\buu \in
\Cl(X)^{\vee}$.  A similar assertion holds for graded right
$A$-modules.
\end{corollary}

\begin{proof}
This condition is clearly necessary.  To see the other direction,
consider the surjective graded morphism defined by the given
generators: $\bigoplus_{i} A(-\bbb_i) \longrightarrow F$.  By
hypothesis, this factors to an epimorphism $\bigoplus_{i} \Dl(-\bbb_i)
\longrightarrow F$ and Corollary~\ref{cor:D_in_M} implies that $F \in
\Ml$.
\end{proof}

\begin{corollary} \label{cor:presentation2}
For each $\cF \in \DMod$, there exist a family $\{\bbb_{i}\}_{i}$ of
elements in $\Cl(X)$ and an epimorphism $\bigoplus_{i} \cD \otimes
\cO(-\bbb_{i}) \longrightarrow \cF$.  The analogous result also holds
for $\cG \in \ModD$.
\end{corollary}

\begin{proof}
Theorem~\ref{thm:Dmodules} implies that there exists $F \in \Ml$ such
that $\cF \cong \tF$.  Now, Corollary~\ref{cor:presentation1} gives an
epimorphism $\bigoplus_{i} \Dl(-\bbb_{i}) \longrightarrow F$.  Taking
the corresponding morphism of sheaves and applying
Corollary~\ref{cor:twisted_D} establishes the claim.
\end{proof}

\begin{remark}
Applying Corollary~\ref{cor:presentation2}, one can construct a
resolution of a $\cD$-module using twisted modules $\cD \otimes
\cO(\bbb_{i})$.  More concretely, given a graded module $F$ satisfying
$\tF \cong \cF$, Gr\"{o}bner basis techniques can be used to construct
a resolution of $F$ by modules $\Dl(\bbb)$ which will then lift to a
resolution of $\cF$.  It would be interesting to investigate the
relationship between these two types of resolutions.
\end{remark}

The last part of this section is devoted to the categorical
equivalence between right and left $\cD$-modules on the toric variety
$X$.  Recall that, for a smooth variety $X$ of dimension $n$, the
sheaf of differential forms of top degree $\Omega^{n}$ has a natural
structure of right $\cD$-module extending the usual $\cO$-module
structure.  Locally, right multiplication of an $n$-form $\omega$ with
a vector field $\nu$ is defined by $\omega \cdot \nu =
-\Lie_{\nu}(\omega)$, where $\Lie_{\nu}(\omega)$ is the Lie derivative
of $\omega$ along $\nu$ (see \cite{borel} Chapter VI.3.4 for details).

The equivalence of categories $\tau_{LR} \colon \DMod \longrightarrow
\ModD$ with inverse $\tau_{RL} \colon \ModD \longrightarrow \DMod$ is
defined as follows: For a left $\cD$-module $\cF$, we have
$\tau_{LR}(\cF) = \cF \otimes \Omega^{n}$ where the right
multiplication with a vector field $\nu$ is $(f \otimes \omega) \cdot
\nu = - \nu(f) \otimes \omega + f \otimes \omega \cdot \nu$.
Similarly, if $\cG$ is a right $\cD$-module, then $\tau_{RL}(\cG) =
{\mathcal H}\!{\it om}_{\cO}(\Omega^{n},\cG)$ and left multiplication
with a vector field $\nu$ is given by $\nu \cdot \psi(\omega) =
\psi(\omega \cdot \nu) - \psi(\omega) \cdot \nu$.

In particular, the left--right equivalence of $A$-modules is given by
the algebra involution $ \tau \colon A \longrightarrow A$, where
$x^{\aa}\partial^{\bb} \mapsto (-\partial)^{\bb} x^{\aa}$.
Specifically, given a graded left $A$-module $F$, we obtain a graded
right $A$-module $F^{\tau}$ which has the same underlying additive
structure and has multiplication defined by $g \cdot f =\tau(f) \cdot
g$ for $g \in A$ and $f \in F$.  Similarly, if $G$ is a graded right
$A$-module, then an analogous procedure yields the left $A$-module
$G^{\tau}$.  It is clear that $(F^{\tau})^{\tau} = F$ and
$(G^{\tau})^{\tau} = G$.  Furthermore, for graded $A$-modules, we have

\begin{proposition} 
There are inverse equivalences of categories
\[
\tau_{LR}^{{\sf mod}} \colon \Ml \longrightarrow \Mr \, , \quad
\tau_{RL}^{{\sf mod}} \colon \Mr \longrightarrow \Ml
\]
given by $\tau_{LR}^{{\sf mod}}(F) = F^{\tau}(-\bee)$ and
$\tau_{RL}^{{\sf mod}}(G) = G^{\tau}(\bee)$ where $\bee \in \Cl(X)$ is
the class of $\ee = \ee_{1} + \dotsb + \ee_{d} \in \ZZ^{d}$.
\end{proposition}

\begin{proof}
We only need to show that the graded components of $\tau_{LR}^{{\sf
mod}}(F)$ and $\tau_{RL}^{{\sf mod}}(G)$ are annihilated by suitable
Euler operators.  However, this follows from the fact that
$\tau(\theta_{\buu}) = - \theta_{\buu} - \langle \buu, \bee \rangle$,
where $\buu \in \Cl(X)$.
\end{proof}

\begin{example}
Since $\tau(\theta_{\buu} + \langle \buu, \baa \rangle) =
-(\theta_{\buu} - \langle \buu, \baa-\bee \rangle)$, there is an
isomorphism $\tau_{LR}^{{\sf mod}}\big(\Dl(\baa) \big) \cong
\Dr(\baa-\bee)$.
\end{example}

We next show that these equivalences of categories are compatible with
the functors in Theorem~\ref{thm:Dmodules}.  We will use the fact that,
for a smooth toric variety $X$, there is a natural isomorphism
$\Omega^{n} \cong \cO(-\bee)$.  In fact, if $0$ denotes the unique
zero dimensional cone in $\Delta$, then this isomorphism identifies
the section $\tfrac{dy_{1}}{y_{1}} \wedge \dotsb \wedge
\tfrac{dy_{n}}{y_{n}}$ with $\frac{1}{x_{1} \dotsm x_{d}}$ on the open
subset $U_{0}$ (see Section~4.3 in \cite{fulton}).

\begin{proposition} \label{prop:compat}
For the pair of functors $(\tau_{LR}^{{\sf mod}},\tau_{LR})$, the
diagrams
\[
\begin{CD}
\Ml @>{\tau_{LR}^{{\sf mod}}}>> \Mr \\
@V{\widetilde{\quad}}VV @V{\widetilde{\quad}}VV \\
\DMod @>{\tau_{LR}}>> \ModD 
\end{CD}
\quad
\begin{CD}
\DMod @>{\tau_{LR}}>> \ModD \\
@V{\Gamma_{L}}VV @V{\Gamma_{R}}VV \\
\Ml @>{\tau_{LR}^{{\sf mod}}}>> \Mr
\end{CD}
\]
are commutative, up to natural isomorphisms.  A similar statement
holds for the pair $(\tau_{RL}^{{\sf mod}}, \tau_{RL})$ .
\end{proposition}

\begin{proof}
Since $\tau_{RL} = (\tau_{LR})^{-1}$ and $\tau_{RL}^{{\sf mod}} =
(\tau_{LR}^{{\sf mod}})^{-1}$, the second assertion is a consequence
of the first.  Because $F^{\tau} \cong F$ as $S$-modules and
$\Omega^{n} \cong \cO(-\bee)$ as $\cO$-modules, there is a natural
isomorphism of $\cO$-modules $\beta_{F} \colon \tau_{LR}(\tF) = \tF
\otimes \Omega^{n} \longrightarrow \widetilde{\tau_{LR}^{{\sf
mod}}(F)} = \widetilde{F^{\tau}(-\bee)}$.  Thus, it suffices to prove
that $\beta_{F}$ is compatible with the right $\cD$-module structures.

By taking a presentation $\bigoplus_{i} \Dl(\bbb_{i}) \longrightarrow
F$, we see that it suffices to establish the claim for $F =
\Dl(\bbb_{i})$.  In this case, the restriction map
\[
H^{0}\big( U, \widetilde{F^{\tau}(-\bee)} \big) \longrightarrow
H^{0} \big( U',\widetilde{F^{\tau}(-\bee)} \big) \, ,
\]
is injective for open subsets $U' \subseteq U \subseteq X$.  Thus, the
claim reduces to showing that $\beta_{F}$ is compatible with the right
$\cD$-module structure on $U_{0}$.  Over $U_{0}$, the map
\[
\beta_{F} \vert_{U_{0}} \colon (F_{x_{1} \dotsm x_{d}})_{\b00}
\otimes_{(S_{x_{1} \dotsm x_{d}})_{\b00}} H^{0}(U_{0},\Omega^{n}) 
\longrightarrow (F_{x_{1} \dotsm x_{d}})_{-\bee}^{\tau} 
\]
is given by $f \otimes \omega \mapsto \tfrac{f}{x_{1} \dotsm x_{d}}$,
where $\omega = \tfrac{dy_{1}}{y_{1}} \wedge \dotsb \wedge
\tfrac{dy_{n}}{y_{n}}$.  Now, it is enough to check that $\beta_{F}
\vert_{U_{0}}$ is compatible with right multiplication with a vector
field $\nu$ over $U_{0}$.  Using the notation from
Lemma~\ref{lem:psi_phi}, we may assume that $\nu = y^{\pp}
\rho(\theta_{i})$, for some $\pp \in \ZZ^{n}$.  We first compute
\[
( f \otimes \omega ) \cdot y^{\pp} \rho(\theta_{i}) = - \big( y^{\pp}
\rho(\theta_{i}) \cdot f \big) \otimes \omega + f \otimes \big( \omega
\cdot y^{\pp} \rho(\theta_{i}) \big) \, .
\]
By definition, we have
\begin{align*}
\begin{split}
( \omega \cdot y^{\pp} \vartheta_{j})&(\eth_{1}, \dotsc, \eth_{n})
= - \big( \Lie_{y^{\pp} \vartheta_{j}} (\omega) \big)(\eth_{1},
\dotsc, \eth_{n}) \\
&= \sum_{i=1}^{n} \omega (\eth_{1}, \dotsc, [y^{\pp} \vartheta_{j},
\eth_{i} ], \dotsc, \eth_{n}) - \Lie_{y^{\pp} \vartheta_{j}}\big( \omega
(\eth_{1}, \dotsc, \eth_{n}) \big) \\
&= - \sum_{i=1}^{n} \left( \delta_{i j} \tfrac{(\pp_{j} +1)
y^{\pp}}{y_{1} \dotsm y_{n}} \right) + \tfrac{y^{\pp}}{y_{1} \dotsm
y_{n}} \\
&= - \tfrac{\pp_{j} y^{\pp}}{y_{1} \dotsm y_{n}} \, ,
\end{split}
\end{align*}
from which we deduce $\omega \cdot y^{\pp} \rho(\theta_{i}) = -
\iota(\pp)_{i} y^{\pp} \cdot \omega$.  Identifying the action of
$y^{\pp} \rho(\theta_{i})$ on $F$ with the action of $x^{\iota(\pp)}
\theta_{i}$, we obtain
\[
( f \otimes \omega ) \cdot y^{\pp} \rho(\theta_{i}) =
\big((-x^{\iota(\pp)} \theta_{i} - \iota(\pp)_{i} x^{\iota(\pp)})
\cdot f \big) \otimes \omega \, .
\]
On the other hand, we have
\begin{align*}
\tfrac{f}{x_{1} \dotsm x_{d}} \cdot x^{\iota(\pp)}\theta_{i} 
&= \tau \big( x^{\iota(\pp)}\theta_{i} \big) \cdot \tfrac{f}{x_{1}
\dotsm x_{d}} \\
&= \big( -x^{\iota(\pp)} \theta_{i} - (\iota(\pp)_{i} +1)x^{\iota(\pp)}
\big) \tfrac{f}{x_{1} \dotsm x_{d}} \\
&= \tfrac{1}{x_{1} \dotsm x_{d}}(- x^{\iota(\pp)} \theta_{i} -
\iota(\pp)_{i} x^{\iota(\pp)}) f \, ,
\end{align*}
and we conclude that $\beta_{F} \big( ( f \otimes \omega ) \cdot \nu
\big) = \beta_{F} ( f \otimes \omega ) \cdot \nu$.

In the second part, we consider $\cF \in \DMod$.  For $F =
\Gamma_{L}(\cF)$, we construct the natural map $\beta_{\cF}' \colon
\tau_{LR}^{{\sf mod}} \big( \Gamma_{L}(\cF) \big) \longrightarrow
\Gamma_{R} \big( \tau_{LR}(\cF) \big)$, by composing the morphisms:
\begin{align} \label{eqn:beta1}
\tau_{LR}^{{\sf mod}} \big( \Gamma_{L}(\cF) \big) = F^{\tau}(-\bee) &
\longrightarrow \Gamma_{R} \big( \widetilde{F^{\tau}(-\bee)} \big) \,
, \\ \label{eqn:beta2}
\beta_{F^{\tau}(-\bee)} \colon \Gamma_{R} \big(
\widetilde{F^{\tau}(-\bee)} \big) & \longrightarrow \Gamma_{R}(\cF
\otimes \Omega^{n}) = \Gamma_{R} \big( \tau_{LR}(\cF) \big) \, .
\end{align}
Since $F$ is $\mathfrak{b}$-saturated, it follows that
$F^{\tau}(-\bee)$ is also $\mathfrak{b}$-saturated and hence
\eqref{eqn:beta1} is an isomorphism.  Moreover, $\Gamma_{R}$ is a
functor and $\beta_{F^{\tau}(-\bee)}$ is an isomorphism which implies
that \eqref{eqn:beta2} is also an isomorphism.
\end{proof}

\begin{example}
The isomorphism of $\cO$-modules $\Omega^{n} \cong \cO(-\bee)$ yields
an isomorphism of $S$-modules $\Gamma_{R}(\Omega^{n}) \cong S(-\bee)$.
Proposition~\ref{prop:compat} shows that this is an isomorphism of
right $A$-modules if $S(-\bee)$ has the right $A$-module structure
given by $S(-\bee) \cong \tfrac{ A(-\bee)}{(\partial_{1}, \dotsc,
\partial_{d}) \cdot A}$.
\end{example}

\section{The Characteristic Variety}

In this section, we use the relationship between $\cD$-modules on $X$
and graded $A$-modules to describe the characteristic varieties.  In
particular, we relate the dimensions of $F$ and $\tF$. For simplicity,
we restrict our attention to left modules.

We start by recalling the quotient construction of $X$; see \cite{cox}
or \cite{musson1}.  Let $T$ be the torus $\Hom\big(\Cl(X),
k^{\ast}\big) \cong (k^{\ast})^{d-n}$.  The group $T$ can be embedded
into $(k^{\ast})^{d}$ by the projection $\ZZ^{d} \longrightarrow
\Cl(X)$.  The diagonal action of $(k^{\ast})^{d}$ on the affine space
$\AA^{d}$ induces an action of $T$ on $\AA^{d}$ such that the open
subset $U = \AA^{d} \setminus \Var(\mathfrak{b})$ is $T$-invariant ---
$\Var(\mathfrak{b})$ denotes the subscheme associated to the ideal
$\mathfrak{b}$.  Since $X$ is smooth (and hence simplicial), there is
a canonical morphism $U \longrightarrow X$ such that $X$ is a
geometric quotient of $U$ with respect to the action of $T$.
Furthermore, we have

\begin{lemma} \label{lem:orbit_dim}
For every $z \in U$, $\Stab_{T}(z) = \{ 1 \}$.  In particular, all the
$T$-orbits in $U$ have dimension $d-n$.
\end{lemma}

\begin{proof}
Consider a point $z$ in $U$ and $t \in T$ satisfying $t \cdot u = u$.
Writing $z = (z_{1}, \dotsc, z_{d})$, we have $t \cdot z = \big(
t(\bee_{1}) z_{1}, \dotsc, t(\bee_{d}) z_{d} \big)$ and we deduce that
$t(\bee_{i}) = 1$, for all $i$ such that $z_{i} \neq 0$.  Because there
is $\sigma \in \Delta$ such that $x^{\hsigma}(z) \neq 0$, we conclude
that $t(\bee_{i}) = 1$ for every $i$ with $\vv_{i} \not\in \sigma$.

On the other hand, $t$ belongs to $\Hom \big( \Cl(X), k^{\ast} \big)$
so we have 
\[
t( \langle \hh, \vv_{1} \rangle \bee_{i} + \dotsb + \langle \hh,
\vv_{d} \rangle \bee_{d}) = t(\bee_{1})^{\langle \hh, \vv_{1} \rangle}
\dotsm t(\bee_{d})^{\langle \hh, \vv_{d} \rangle} = 1 \, ,
\]
for every $\hh \in N^{\vee}$.  It follows that $\prod_{\vv_{i} \in
\sigma} t(\bee_i)^{\langle \hh, \vv_{i} \rangle} = 1$.  Because this
holds for each $\hh \in N^{\vee}$ and the $\vv_{i}$ form part of a
basis of $N$, we also conclude that $t(\bee_{i}) = 1$ when $v_{i} \in
\sigma$ and therefore $t = 1$.
\end{proof} 

We next identify $\AA^{d} \times \AA^{d}$ with the cotangent bundle of
$\AA^{d}$ and consider the natural $T$-action on it.  Let $S' =
k[x_{1}, \dotsc, x_{d}, \xi_{1}, \dotsc, \xi_{d}]$, with the
$\Cl(X)$-grading given by $\deg(x_{i}) = - \deg(\xi_i) = \bee_{i}$, be
the coordinate ring of $\AA^{d} \times \AA^{d}$.  Since the action of
$T$ on $\AA^{d}$ is linear, it follows that the action of $T$ on
$\AA^{d} \times \AA^{d}$ is given by $t \cdot (z_{1}, z_{2}) = (t
\cdot z_{1}, t^{-1} \cdot z_{2})$.  It is clear that $V = U \times
\AA^{d} \subset \AA^{d} \times \AA^{d}$ is invariant under the action
of $T$.  We construct its quotient as follows.

\begin{proposition} \label{prop:quotient}
There is a morphism $\pi \colon V \longrightarrow X'$ such that $X'$
is the geometric quotient of $V$ by the action of $T$.  In addition,
for every $z \in V$, we have $\Stab_{T}(z) = \{ 1 \}$ implying that
all the $T$-orbits in $V$ have dimension $d-n$.
\end{proposition}

\begin{proof}
The first step is to construct the morphism $\pi \colon V
\longrightarrow X'$ as a categorical quotient --- this is a local
problem.  For every $\sigma \in \Delta$, let $V_{\sigma} \subseteq V$
be the open subset defined by the non-vanishing of $x^{\hsigma}$.  In
other words, we have $V_{\sigma} = \big( U \setminus \Var(x^{\hsigma})
\big) \times \AA^{d} \subseteq V$, which is clearly $T$-invariant.
Thus, the categorical quotient is locally $V_{\sigma} = \Spec\big(
S'[(x^{\hsigma})^{-1}]^{T} \big)$.  Since $t \cdot x^{\aa} \xi^{\bb} =
t(\baa-\bbb) x^{\aa}\xi^{\bb}$, for every $t \in T$ and $\aa$, $\bb
\in \ZZ^{d}$, we have $S'[(x^{\hsigma})^{-1}]^{T} =
S'[(x^{\hsigma})^{-1}]_{\b00}$.  Now, if $\sigma_{0}$ is a face of
$\sigma$ such that $\sigma_{0} = \sigma \cap \hh^{\perp}$ for $\hh \in
N^{\vee} \cap \sigma^{\vee}$, then we set $\cc = \langle \hh, \vv_{1}
\rangle \ee_{1} + \dotsb + \langle \hh, \vv_{d} \rangle \ee_{d} \in
\ZZ^{d}$.  It follows that $S'[(x^{\hsigma_{0}})^{-1}]_{\b00} = \big(
S'[(x^{\hsigma})^{-1}]_{\b00} \big)_{x^{\cc}}$ which provides an open
immersion $V_{\sigma_{0}} \hookrightarrow V_{\sigma}$.  Thus, we
obtain morphisms $\pi_{\sigma} \colon V_{\sigma_{0}} \longrightarrow
V_{\sigma}$ which glue together to give the categorical quotient $\pi
\colon V \longrightarrow X'$.

In the second step, we establish that $\pi \colon V \longrightarrow
X'$ is in fact a geometric quotient.  By Amplification~1.3 in
\cite{fmu}, it suffices to show that every $T$-orbit in $V$ is closed.
Consider $z = (z_{1}, z_{2}) \in V$.  Since $U \longrightarrow X$ is a
geometric quotient, the projection $V \longrightarrow U$ induces a
morphism $\chi \colon \overline{Tz} \longrightarrow \overline{Tz_{1}}
= Tz_{1}$.  By Lemma~\ref{lem:orbit_dim}, the morphism $\gamma \colon
T \longrightarrow Tz_{1}$ given by $\gamma(t) = tz_{1}$ is bijective.
Because the characteristic of the ground field $k$ is zero and both
$T$ and $Tz_{1}$ are smooth, the map $\gamma$ is an isomorphism.  Let
$\gamma' \colon T \longrightarrow Tz \subseteq V$ be defined by
$\gamma'(t) = tz$.  Hence, the map $\gamma' \circ \gamma^{-1} \circ
\chi \colon \overline{Tz} \longrightarrow Tz \subseteq \overline{Tz}$
is the identity on $Tz$.  It follows that $\gamma' \circ \gamma^{-1}
\circ \chi$ is the identity map on $\overline{Tz}$ and we conclude
that $Tz = \overline{Tz}$.

Since the projection from $V$ onto $U$ is $T$-equivariant, the second
assertion follows from Lemma~\ref{lem:orbit_dim}.
\end{proof}

Before discussing characteristic varieties, we review some properties
of the order filtration.  Recall that the sheaf $\cD$ is naturally
filtered by the order of the differential operators.  In particular,
this makes $H^{0}(U_{\sigma}, \cD)$ into a filtered ring.  We can also
filter the ring $A$ by the order of the differential operators (in
fact, $S' = \gr(A)$) and this induces a filtration on the quotient
$\frac{(A_{x^{\hsigma}})_{\b00}}{(A_{x^{\hsigma}})_{\b00} \cdot
(\theta_{\buu} : \buu \in \Cl(X)^{\vee})}$.  As Musson observed, we
have

\begin{lemma}[Musson] \label{lem:filtered_map}
For every $\sigma \in \Delta$, the isomorphism
$\bar{\varphi}_{\sigma}$ $[$see equation~\eqref{eqn:varphi_map}$]$
preserves the filtrations induced by the order of differential
operators.
\end{lemma}

\begin{proof}
See Section~4 in \cite{musson1}.
\end{proof}

A filtration of an $A$-module $F$ is called good if the associated
graded module $\gr(F)$ is finitely generated over $S'$.  Every
finitely generated $A$-module has a good filtration and, conversely,
any module with a good filtration is necessarily finitely generated
over $A$.  For a good filtration of $F$, we define the characteristic
ideal $\mathfrak{i}(F)$ to be the radical of $\Ann_{S'}\big( \gr(F)
\big)$.  Since any two good filtrations are equivalent, the
characteristic ideal $\mathfrak{i}(F)$ is independent of the choice of
good filtration.  The characteristic variety of $F$ is $\CharVar(F) =
\Var\big( \mathfrak{i}(F) \big) \subseteq \AA^{d} \times \AA^{d}$.
Analogously, for a $\cD$-module $\cF$ with a good filtration, we
define the characteristic variety $\CharVar(\cF)$ to be the support of
the associated graded sheaf $\gr(\cF)$.

We first describe the characteristic variety associated to the graded
left $A$-modules $\Dl(\bbb)$.  Let $p_{\buu} = \langle \buu, \bee_{1}
\rangle x_{1}\xi_{1} + \dotsb + \langle \buu, \bee_{d} \rangle
x_{d}\xi_{d} \in S'$, for all $\buu \in \Cl(X)^{\vee}$.  We consider
the ideal $\mathfrak{p} = (p_{\buu} : \buu \in \Cl(X)^{\vee})$ and the
corresponding variety $Z = \Var(\mathfrak{p}) \subseteq \AA^{d} \times
\AA^{d}$.  It is clear that $Z$ is invariant under the $T$-action.

\begin{proposition} \label{prop:char_D}
The variety $Z$ is a normal complete intersection of dimension $d+n$.
Moreover, $Z$ is equal to the characteristic variety
$\CharVar\big(\Dl(\bbb)\big)$, for every $\bbb\in\Cl(X)$.
\end{proposition}

\begin{proof}
By choosing a basis $\buu_{1}, \dotsc, \buu_{d-n}$ for
$\Cl(X)^{\vee}$, we can write $\mathfrak{p} = (p_{\buu_{i}} : 1 \leq i
\leq d-n)$.  For each $\buu_{i}$, we pick a representative $\uu_{i}
\in (\ZZ^{d})^{\vee}$ for $\buu_{i}$ such that $\uu_{1}, \dotsc,
\uu_{d-n}$ are linearly independent.  We then enlarge this collection to
obtain a basis $\uu_{1}, \dotsc, \uu_{d}$ for $(\ZZ^{d})^{\vee}$.
Setting $q_{\uu} = \langle \uu, \ee_{1} \rangle x_{1}\xi_{1} + \dotsb
+ \langle \uu, \ee_{d} \rangle x_{d}\xi_{d} \in S'$, for all $\uu \in
(\ZZ^{d})^{\vee}$, it follows that the ideal $(q_{\uu_{i}} : 1 \leq i
\leq d)$ equals $(x_{i}\xi_{i} : 1 \leq i \leq d)$, which has
dimension $d$.  Hence, the $q_{\uu_{i}}$ and $p_{\buu_{i}}$ form a
regular sequence and we deduce $\dim(Z) = d+n$.

To prove that $Z$ is normal, we apply Serre's criterion.  Because
being a complete intersection implies the $(S2)$ condition, it
suffices to show that $Z$ satisfies condition $(R1)$, which we check
by using the Jacobian criterion.  The Jacobian matrix $\Jac(x, \xi)$
of $(p_{\buu_{1}}, \dotsc, p_{\buu_{d-n}})$ is given by
\[
\begin{pmatrix}
\langle \buu_{1}, \bee_{1} \rangle \xi_{1} & \dotsb & 
\langle \buu_{1}, \bee_{d} \rangle \xi_{d} & 
\langle \buu_{1}, \bee_{1} \rangle x_{1} & \dotsb & 
\langle \buu_{1}, \bee_{d} \rangle x_{d} \\
\vdots & \ddots & \vdots & \vdots & \ddots & \vdots \\
\langle \buu_{d-n}, \bee_{1} \rangle \xi_{1} & \dotsb & 
\langle \buu_{d-n}, \bee_{d} \rangle \xi_{d} & 
\langle \buu_{d-n}, \bee_{1} \rangle x_{1} & \dotsb & 
\langle \buu_{d-n}, \bee_{d} \rangle x_{d} 
\end{pmatrix} \, .
\]
Observe that, for $1 \leq i \leq d$, the restriction $\bigoplus_{j
\neq i} \ZZ \ee_{j} \longrightarrow \Cl(X)$ is surjective.  Indeed, if
$\sigma_{i}$ is the cone generated by $\vv_{i}$, then every element in
$\Cl(X)$ can be represented by a divisor whose support does not
intersect $U_{\sigma_{i}}$.  We deduce that if the rank of $\Jac(x, \xi)$
is strictly less than $d-n$, then at least two of the pairs of
coordinates $(x_{1},\xi_{1}), \dotsc, (x_{d}, \xi_{d})$ are zero.  By
cutting with $n$ extra quadrics $q_{\uu_{d-n+1}}, \dotsc, q_{\uu_{d}}$,
we see that the codimension of the singular locus of $Z$ is at least
two.  Therefore, the variety $Z$ is normal.  Moreover $Z$ is a cone
which implies it is connected and hence integral.

Recall that, for $\bbb \in \Cl(X)^{\vee}$, we have $ \Dl(\bbb) =
\tfrac{A(\bbb)}{(\theta_{\buu_{i}} + \langle \buu_{i}, \bbb \rangle :
\text{$1\leq i\leq d-n)$}}$, and for the order filtration the initial
term (or principal symbol) of the above elements is
$\IN(\theta_{\buu_{i}} + \langle \buu_{i}, \bbb \rangle) =
\theta_{\buu_{i}}$.  Since the $\theta_{\buu_{i}}$ for $1 \leq i \leq
d-n$ form a regular sequence in $S' = \gr(A)$, it follows that
$\gr\big(\Dl(\bbb)\big) = S'/\mathfrak{p}$.  On the other hand, we
have already seen that $\mathfrak{p}$ is reduced, so we have
$\CharVar\big(\Dl(\bbb)\big) = \Var(\mathfrak{p}) = Z$.
\end{proof}

\begin{corollary} 
If $F \in \Mlf$, then we have $\CharVar(F) \subseteq Z$.
\end{corollary}

\begin{proof}
By Corollary~\ref{cor:presentation1}, $F$ is a quotient of
$\bigoplus_{i=1}^{r} \Dl(\bbb_{i})$, for some $r$ and some $\bbb_{i}$.
It follows that $\CharVar(F) \subseteq \bigcup_{i=1}^{r}
\CharVar\big(\Dl(\bbb_{i})\big) = Z$.
\end{proof} 

We next relate the cotangent bundle of $X$ to the variety
$Z$. Consider the following diagram:
\[
\begin{CD}
Z \cap V @>>> V @>>> U \\
@VVV @V{\pi}VV @VVV \\
\pi(Z \cap V) @>>> X' @>>> X 
\end{CD}
\]
where $X' \longrightarrow X$ arises from the universal property of the
categorical quotient.

\begin{proposition} \label{prop:cotangent}
There is a canonical isomorphism of varieties over $X$ between $\zeta
\colon \pi(Z \cap V) \longrightarrow X$ and the cotangent bundle
$T^{*}X$ over $X$.
\end{proposition}

\begin{proof}
Since $T^{*}X$ is naturally isomorphic to $\CharVar(\cD)$, we see that
$T^{*}X$ is isomorphic to $\Spec\big(\gr H^{0}(U_{\sigma},\cD) \big)$
over $U_{\sigma}$.  On the other hand, from the local description of
$X'$, we know that the inverse image of $U_{\sigma}$ is $\Spec
(S'_{x^{\hsigma}})_{\b00}$ and therefore
\[
\zeta^{-1}(U_{\sigma}) =
\Spec\left(\frac{(S'_{x^{\hsigma}})_{\b00}}{(p_{\buu} : \buu \in
\Cl(X)^{\vee})} \right) \, .
\]
By Lemma~\ref{lem:filtered_map}, we have an isomorphism of filtered
rings:
\[
\bar{\varphi}_{\sigma} \colon
\frac{(A_{x^{\hsigma}})_{\b00}}{(A_{x^{\hsigma}})_{\b00}
\cdot(\theta_{\buu} : \buu \in \Cl(X)^{\vee})} \longrightarrow
H^{0}(U_{\sigma},\cD) \, .
\]
Notice that the graded ring associated to left hand side is
$\frac{(S'_{x^{\hsigma}})_{\b00}} {(p_{\buu} : \buu \in
\Cl(X)^{\vee})}$.  Indeed, following the proof of
Proposition~\ref{prop:char_D}, the initial terms of $\theta_{\buu_1},
\dotsc, \theta_{\buu_ {d-n}}$ are equal to $p_{\buu_1}, \dotsc,
p_{\buu_{d-n}}$ and form a regular sequence in
$(S'_{x^{\hsigma}})_{\b00}$.  Therefore, by passing to the associated
graded rings, $\bar{\varphi}_{\sigma}$ induces the required isomorphism.
Because the $\bar{\varphi}_{\sigma}$ are compatible with restriction, these
local isomorphisms glue together to give the required isomorphism
\end{proof}

We now present the main result in this section.

\begin{theorem} \label{thm:char=char}
If $F \in \Mlf$, then the characteristic variety of $F$ is
$T$-invariant and $\pi \big(\CharVar(F) \setminus \Var(\mathfrak{b})
\times \AA^{d} \big) = \CharVar(\tF)$.
\end{theorem}

\begin{proof}
Since $F \in \Mlf$, we may choose a finite set $f_{1}, \dotsc, f_{r}$
of homogeneous generators for $F$.  By using these homogeneous
elements to define a good filtration of $F$, it follows that $\gr(F)$
is a graded finitely generated $S'$-module.  Therefore, both
$\mathfrak{j}(F) = \Ann_{S'}\big( \gr(F) \big)$ and its radical
$\mathfrak{i} = \sqrt{\mathfrak{j}(F)}$ are graded ideals of $S'$.
Recall that for every $t \in T$, we have $t \cdot x^{\aa}\xi^{\bb} =
t(\baa-\bbb)x^{\aa}\xi^{\bb}$.  We deduce that every subscheme defined
by a graded ideal is $T$-invariant; in particular, $\CharVar(F) =
\Var(\mathfrak{i})$ is $T$-invariant.

To prove the second assertion, we argue locally and use the
identification in Proposition~\ref{prop:cotangent}.  Over the open
subset $U_{\sigma}$, the ideal defining $\CharVar(F) \setminus
\Var(\mathfrak{b}) \times \AA^{d} \subseteq \Spec (S'_{x^{\hsigma}})$
is $\mathfrak{j} \cdot S'_{x^{\hsigma}}$.  On the other hand, the
characteristic variety of $\tF$ over $U_{\sigma}$ can be computed as
follows: If $H^{0}(U_{\sigma}, \tF) = (F_{x^{\hsigma}})_{\b00}$ has
the good filtration induced by the images of $f_{1}, \dotsc, f_{r}$,
then we obtain $\gr\big(H^{0}(U_{\sigma},\tF)\big) \cong
(\gr(F)_{x^{\hsigma}})_{\b00}$.  To see that the annihilator of
$\gr\big(H^{0}(U_{\sigma},\tF)\big)$ in $(S'_{x^{\hsigma}})_{\b00}$ is
$(\mathfrak{j}_{x^{\hsigma}})_{\b00}$, it is enough to observe that
$\gr(F)_{x^{\hsigma}}$ can be generated by elements of degree zero.
Since $\sqrt{(\mathfrak{i}(F)_{x^{\hsigma}})_{\b00}} =
(\mathfrak{j}_{x^{\hsigma}})_{\b00}$, we may identify $\pi(\CharVar(F)
\cap V_{\sigma})$ with $\CharVar(\tF \vert_{U_{\sigma}})$.  Because
these identifications are compatible with restriction, they glue
together to give the required isomorphism.
\end{proof} 

As a corollary, we obtain 

\begin{proof}[Proof of Theorem~\ref{thm:maincharvar}]
Follows immediately from Theorem~\ref{thm:char=char}.
\end{proof}

We end by relating the dimension of the $A$-module $F$ and its
associated $\cD$-module $\tF$.  By definition, the local dimension of
a $\cD$-module $\cF$ at a point $p \in X$ is equal to the Krull
dimension of the associated graded module of $\cF_p$ with respect to a
good filtration respecting the order filtration of $\cD_{p}$.  It is
also equal to the dimension of the characteristic variety of $\cF_p$.
The dimension of $F$ is by definition the maximum of the local
dimensions and is equivalently the dimension of $\CharVar(F)$.

By Theorem~\ref{thm:char=char}, the dimension of $\tF$ is equal to the
maximum of the local dimensions of $F$ over the open set $\AA^{d}
\setminus \Var(\mathfrak{b})$, minus the dimension $d-n$ of the orbits
under the group action.  Before showing that when $F$ has no
$\mathfrak{b}$-torsion, we can express this dimension in terms of
$\dim(F)$, we collect two lemmas.

\begin{lemma} \label{lem:local_dim}
If $F$ is a finitely generated left $A$-module, $f$ is an element of
$S$, and $F'$ is a finitely generated $A$-submodule of $F[f^{-1}]$,
then $\dim(F') \leq \dim(F)$.
\end{lemma}

\begin{proof}
This claim follows immediately from standard results about
Gelfand-Kirillov dimension; see Propositions~8.3.2~(i) and
8.3.14~(iii) in \cite{McR}.
\end{proof}


\begin{proposition} \label{prop:torsionfree_dim}
Let $F \in \Mlf$ and recall that $\mathfrak{b}$ is the irrelevant
ideal in $S$.  If $F$ has no $\mathfrak{b}$-torsion, then we have 
\[
\dim(F) =\max\{ \dim(F_{p}) : \text{$p \in \AA^{d} \setminus
\Var(\mathfrak{b})$} \} \, .
\]
\end{proposition}

\begin{proof}
Suppose otherwise; then we have $\dim(F) > \dim(F_{z})$ for all $z \in
\AA^{d} \setminus \Var(\mathfrak{b})$.  Let $F'$ be the maximum
submodule of $F$ of dimension strictly less than $\dim(F)$.  In other
words, $F'$ is the submodule consisting of all $f \in F$ such that $A
\cdot f$ has dimension strictly less than $\dim(F)$.  Since $F'$ is a
submodule, there is a short exact sequence
\[
0 \longrightarrow F' \longrightarrow F \longrightarrow \tfrac{F}{F'}
\longrightarrow 0 \, .
\]
By construction, $F / F'$ has no nonzero submodules of dimension
strictly less than $\dim(F)$ and, hence, the irreducible components of
$\CharVar(F / F')$ have dimension at least $\dim(F)$; see
\cite{smith}.  By hypothesis, the irreducible components of
$\CharVar(F)$ of dimension $\dim(F)$ are contained inside
$\zeta^{-1}(\Var(\mathfrak{b}))$ where $\zeta : T^{\ast}X
\longrightarrow X$.  Since $\CharVar(F) = \CharVar(F') \cup
\CharVar(F/F')$, it follows that the characteristic variety of $F/F'$
is contained inside $\zeta^{-1}(\Var(\mathfrak{b}))$.  Moreover, the
support of an $A$-module equals the projection of its characteristic
variety (see \cite{gm}) which implies that $F/F'$ is supported on
$\Var(\mathfrak{b})$.  Taking the long exact sequence in local
cohomology, we have
\[
0 \longrightarrow H^{0}_{\mathfrak{b}}(F') \longrightarrow
H^{0}_{\mathfrak{b}}(F) \longrightarrow
H^{0}_{\mathfrak{b}}\left(\tfrac{F}{F'}\right) \longrightarrow
H^{1}_{\mathfrak{b}}(F') \longrightarrow \dotsb \, .
\]
By assumption $F$ has no $\mathfrak{b}$-torsion, so we have
$H^{0}_{\mathfrak{b}}(F') = H^{0}_{\mathfrak{b}}(F) = 0$.  Because
$F/F'$ is supported on $\Var(\mathfrak{b})$, we have
$H^{0}_{\mathfrak{b}}(F/F') = F/F'$.  Choosing a set of generators
$\mathfrak{b} = (s_{1}, \dotsc, s_{r})$, the long exact sequence
induces
\[
0 \longrightarrow \tfrac{F}{F'} \longrightarrow
\tfrac{\bigoplus_{i=1}^{r} F'[s_{i}^{-1}]}{F'} \longrightarrow \dotsb
\, .
\]
Hence, $F/F'$ is a finitely generated $A$-subquotient of
$\bigoplus_{i=1}^{r} F'[s_{i}^{-1}]$.  Lemma~\ref{lem:local_dim} then
implies that $\dim(F/F') \leq \dim(F') < \dim(F/F')$ which is a
contradiction.
\end{proof}

Finally, we have

\begin{theorem} \label{thm:dim}
If $F \in \Mlf$ has no $\mathfrak{b}$-torsion, then we have $\dim(\tF)
= \dim(F) - d + n$.
\end{theorem}

\begin{proof}
Applying Proposition~\ref{prop:torsionfree_dim}, we see that $\dim(F)$
is the maximum of the local dimensions of $F$ over $\AA^{d} \setminus
\Var(\mathfrak{b})$.  Hence, the claim follows from
Proposition~\ref{prop:quotient}.
\end{proof}

A coherent $\cD$-module $\cF$ is holonomic if $\dim(\cF) = \dim X$.

\begin{corollary}
If $F \in \Mlf$ is holonomic, then $\tF$ is holonomic.  Furthermore,
every holonomic $\cD$-module is of the form $\tF$ for some holonomic
$F \in \Mlf$.
\end{corollary}

\begin{proof}
Follows immediately from Theorem~\ref{thm:dim} and
Remark~\ref{rem:rep}.
\end{proof}


\newcommand{\etalchar}[1]{$^{#1}$}

\end{document}